\documentclass[preprint,groupedaddress,amsmath,amssymb,showkeys]{revtex4}
\usepackage{dcolumn}% Align table columns on decimal point
\usepackage{bm}% bold math

\usepackage{amsmath}
\usepackage{amssymb}
\usepackage{latexsym}
\usepackage[dvips]{graphicx}
\usepackage{afterpage}
\usepackage{float}
\newcommand{\m}{\mbox}
\newcommand{\bd}{\boldmath}
\newcommand{\bo}[1]{\mbox{\boldmath $#1$}}

\begin{document}

%\preprint{}

\title{Helical buckling of a whirling conducting rod in a uniform magnetic field}


\author{J.~Valverde}
%\email{}
%\thanks{}
\affiliation{Department of Civil and Environmental Engineering, University
of California, Berkeley, USA}
\author{G.H.M.~van~der~Heijden}
\email{g.heijden@ucl.ac.uk}
%used with the \author command). \
\affiliation{Centre for Nonlinear Dynamics, University College London,\\
Gower Street, London WC1E 6BT, UK}

\date{\today}

\begin{abstract}
We study the effect of a magnetic field on the behaviour of a
conducting elastic rod subject to a novel set of boundary conditions
that, in the case of a transversely isotropic rod, give rise to exact
helical post-buckling solutions. The equations used are the geometrically
exact Kirchhoff equations and both static (buckling) and dynamic (whirling)
instability are considered. Critical loads are obtained explicitly and are
given by a surprisingly simple formula. By solving the linearised equations
about the (quasi-)stationary solutions we also find secondary instabilities
described by \mbox{(Hamiltonian-)Hopf} bifurcations, the usual signature of
incipient `breathing' modes. The boundary conditions can also be used to
generate and study helical solutions through traditional non-magnetic
buckling due to compression, twist or whirl.
\end{abstract}


\maketitle

%\section{}
%\subsection{}
%\subsubsection{}

\section{Introduction}

\begin{figure}
\begin{center}
\includegraphics[width=0.4\linewidth]{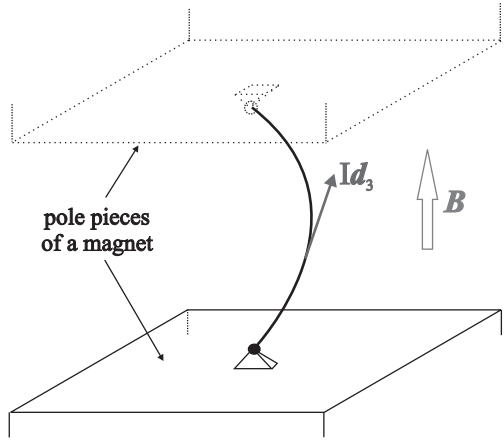}
\end{center}
\caption{Experimental setup for a conducting wire.}
\label{fig:hel:exp-setup}
\end{figure}

A straight current-carrying wire held in tension between pole faces
of a magnet is well known to buckle into a (roughly) helical
configuration at a critical current (see
Fig.~\ref{fig:hel:exp-setup}). A photograph of this phenomenon is
shown in Section 10.4.3 of \cite{woodson}, where
a linear stability analysis is carried out for a simple string
model. (A string is here meant to be a perfectly flexible elastic
wire.) The problem was studied by Wolfe \cite{wolfe1} by means of a
rigorous bifurcation analysis for a (nonlinearly elastic) string
suspended between fixed supports and placed in a uniform magnetic
field directed parallel to the undeformed wire. He found that an
infinite number of solution branches bifurcate from the trivial
straight solution, much like in the Euler elastica under compressive
load. In this case the non-trivial solutions are exact helices. That
this should be so, is easily explained by the fact that the
(Lorentz) body force is everywhere normal to the deformed
configuration and hence the wire necessarily in a uniform state of
tension. Some (statics) stability results (i.e., minimisation of
the potential energy) were also obtained, indicating that the
first branch of solutions is stable while the others are unstable.

In a subsequent paper Wolfe \cite{wolfe2} extends the analysis to a
uniformly rotating (whirling) string and shows again the existence
of bifurcating branches of whirling non-trivial solutions. Due to
centrifugal effects no closed-form solutions could be obtained in
this case. This result was further extended by Healey \cite{healey}
using equivariant bifurcation theory in order to deal with the
symmetries of the problem.

Wolfe also considered a conducting {\it rod} in a uniform magnetic
field \cite{wolfe3}. In addition to extension a rod can undergo
flexure, torsion and shear, and for the case of welded boundary
conditions it was found that in certain cases bifurcation occurs,
with the usual infinity of non-trivial equilibrium states. All the
works cited above were content with showing the existence of
bifurcating solutions and did not study their post-buckling
behaviour.

In this paper we consider the post-buckling behaviour of a
conducting rod. Wolfe considered welded boundary conditions in
\cite{wolfe3} and encountered degeneracies (even-dimensional
eigenspaces) because of rotational symmetry of the problem.
In previous work \cite{clamped-JNLS} we showed that further
complications occur and that magnetic buckling of a welded
transversely isotropic rod (i.e., a rod with unequal bending
stiffnesses about the two principal axes of its cross-section)
is described by a remarkably degenerate pitchfork bifurcation.
Wolfe also reported numerical evidence of helical post-buckling
solutions. However, exact helical solutions cannot be supported by
(coaxial) welded boundary conditions. Here we formulate a novel set
of what we call `coat hanger' boundary conditions that do support
(i.e., are compatible with) exact helical solutions, and show that
subject to these boundary conditions an isotropic rod does indeed
buckle (exclusively) into a helix, or more precisely, that there is
an infinite series of helical modes bifurcating at increasing load,
each successive mode having one more (half) helical turn. All helical
solutions can be obtained explicitly and it is found that the
pitchfork bifurcations for these coat hanger boundary conditions
are non-degenerate and that the critical loads are given by a
remarkably simple formula. Unlike in string buckling a rod does not
require a tensile force in the trivial state, but we allow for such
an applied force as well. The pertinent dimensionless parameter that
governs buckling measures the product of current and magnetic field
against the bending force.

We also study steady whirling solutions for which we introduce a
rotating coordinate system. This extends Wolfe's analysis of whirling
strings to whirling rods. An interesting feature of helical solutions
is that since all points on a helix have equal distance to the whirling
axis, and are therefore equally affected by centrifugal forces, solutions
remain helical when spun. We perform a stability analysis by computing
eigenvalues of the linearised boudary-value problem about a (quasi-stationary)
whirling solution. For this we use a continuation (or homotopy) approach
that takes advantage of the fact that exact expressions for the (imaginary)
eigenvalues can be obtained in an appropriate limit (no spin, no magnetic
field). The eigenvalues in this limit are then traced as system parameters
are varied.

Whirl tends to destabilise the helical solutions, but stable solutions can
be obtained by adding the effect of internal viscoelastic damping. We find
Hopf bifurcations on the first bifurcating branch where a stable whirling
solution becomes unstable under an increase of the angular velocity. We
also briefly consider anisotropic rods. Critical loads can still be obtained
analytically, but these rods buckle into coiled but non-helical solutions.
Secondary instabilities are found due to Hamiltonian-Hopf bifurcations, a
common signature of `breathing' or `flutter' instabilities in mechanical
systems.

Helical solutions are widely studied in a whole range of applications.
Often these solutions are thought to arise through buckling of a straight
rod under the action of end loads. However, as commented above, exact
helical solutions are not supported by the usual set of boundary conditions.
Consequently, boundary conditions are often not mentioned, or the rod is
implicitly assumed to be infinitely long in order to prevent end effects
\cite{goriely}. One of the contributions of this paper is to present and
highlight boundary conditions for an elastic rod that do support helical
solutions. One could apply these experimentally if one was interested in
generating or studying helical solutions in a finite-length rod.

The paper is organised as follows. In Section 2 we present the rod mechanics
formulation, in which the magnetic field enters the force balance equation
through the Lorentz body force. The coat hanger boundary conditions are
introduced together with a sketch of a device that can be constructed to
realise these boundary conditions in a testing rig. For the study of whirling
solutions the equilibrium equations are transformed to a coordinate system
rotating at constant angular velocity. We use numerical bifurcation and
continuation methods to find the buckling loads and to compute post-buckling
solution paths, both for the statics and dynamics case. After the
nondimensionalisation in Section 3, the linearisation is presented in Section
4. Section 5 first presents analytical stability results for the statics case
and then introduces our continuation approach to numerical stability analysis
of the full system. In Section 6 results are presented in the form of
bifurcation diagrams and curves of Hopf bifurcations in appropriate parameter
planes. Conclusions are drawn in Section 7 and the study closes with two
Appendices: one giving details about the linearised system of equations and
one deriving exact buckling results for helical solutions. The latter is
complementary to the bifurcation analysis (not assuming any shape) in
Section 5; together these analyses give a complete picture of helical
magnetic buckling.

\section{The rod mechanics model}

We describe the elastic behaviour of a conducting cable by the
Kirchhoff equations for the dynamics of thin rods. The rod is
assumed to be uniform, inextensible, unshearable and intrinsically
straight and prismatic. The assumptions of inextensibility and
unshearability are appropriate for thin rods with relatively low
external (here electrodynamic) forces. For the background of the
Kirchhoff equations the reader is referred to \cite{antman,coleman_et_al}.
These equations were also used in \cite{clamped-JNLS} and
\cite{valverde-cosserat} to analyse the dynamics of a spinning tether.

\begin{figure}
\begin{center}
\includegraphics[width=0.5\linewidth]{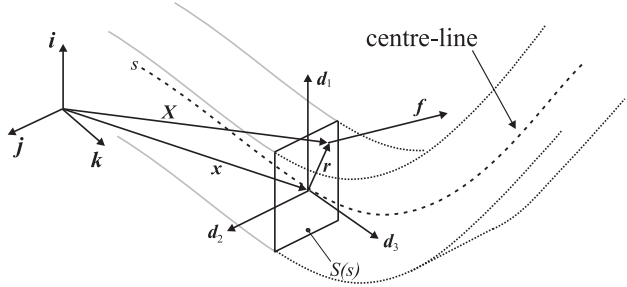}
\end{center}
\caption{Cosserat model of a rod.} \label{fig:cosserat}
\end{figure}

Let $\bo x$ denote the position of the rod's centreline and let
$\{\bo d_1,\bo d_2,\bo d_3\}$ be a right-handed orthonormal frame of
directors (the Cosserat triad) defined at each point along the
centreline. Since the centreline is assumed to be inextensible we
can take $\bo d_3$ in the direction of the local tangent:
\begin{equation}
\bo x'(s,t)=\bo d_3(s,t), \label{eqtangent1}
\end{equation}

\noindent where the prime denotes differentiation with respect to
arclength $s$ measured along the centreline, and $t$ is time. The
directors $\bo d_1$ and $\bo d_2$ will be taken to point along the
principal bending axes of the cross-section (see
Fig.~\ref{fig:cosserat}). The unstressed rod is taken to lie along
the basis vector $\bo k$ of a fixed inertial frame $\{\bo i,\bo
j,\bo k\}$.

Looking at Fig.~\ref{fig:cosserat} we note that the position
vector of an arbitrary point of the rod can be expressed as
\begin{equation}
\begin{split}
\bo X (s,\xi_1,\xi_2,t)=&\bo x (s,t) +\xi_1\bo d _1(s,t) +\xi_2\bo d
_2(s,t)\\
=&\bo x (s,t) +\bo r (s,\xi_1,\xi_2,t), \label{eq:kinem}
\end{split}
\end{equation}

\noindent where $(\xi_1,\xi_2)$ are the components of $\bo r$ in the
cross-section relative to $\{\m{\bd$d$}_1(s),\m{\bd$d$}_2(s)\}$. The
rod is thus viewed as a set of infinitesimal slices centred at all
$s$. A one-dimensional description will be obtained by averaging of
forces and moments over each cross-section. The internal traction,
which is the projection of the stress tensor onto the
cross-sectional plane, is given by a force which we denote by $\bo f
=\bo f (s,\xi_1,\xi_2,t)$ (see Fig.~\ref{fig:cosserat}). The
resultant elastic force exerted in a section $S(s)$ is given by
\begin{equation}
\bo F (s,t)=\int_{S(s)} \bo f (s,\xi_1,\xi_2,t)\,\textrm{d}S,
\end{equation}
\noindent where $\textrm{d}S$ is an infinitesimal area element. This
force can be expressed in the director basis as $\bo F =
\sum_{i=1}^{3} F_{i}\,\m{\bd$d$}_i$. The resultant moment in the
section $S(s)$ is given by
\begin{equation}
\bo M (s,t)=\int_{S(s)} \bo r (s,\xi_1,\xi_2,t) \times  \bo f
(s,\xi_1,\xi_2,t)\,\textrm{d}S, \label{moments}
\end{equation}
\noindent and will be expressed as $\bo M = \sum_{i=1}^{3}
M_{i}\m{\bd$d$}_i$.

The rod is assumed to carry an electric current for which we can
write
\begin{equation}
\bo I = I \bo x' = I \bo d_3.
\end{equation}
Here we have assumed the current to have the same direction as the
rod, which is consistent with a one-dimensional rod theory. It
amounts to the assumption that the cross-section of the conducting
wire is small enough to make currents within the cross-section (eddy
currents) induced by the motion negligible. The current $\bo I$
interacts with the magnetic field $\bo B_0$ to generate a (Lorentz)
body force given by
\begin{equation}
\bo F_L= I\bo d_3 \times \bo B_0.
\end{equation}
Following \cite{wolfe1} we assume the magnetic field to be uniform
and directed along the unstressed rod, i.e.,
\begin{equation}
\bo B_0= B_0 \bo k.
\end{equation}

The balancing of forces and moments across an infinitesimal rod
element then yields the following set of partial differential equations
\cite{antman,coleman_et_al}:
\begin{equation}
\bo F' + IB_0 \bo d_3 \times \bo k = \rho A  \ddot{\bo x} ,
\label{eq:lmomentum2}
\end{equation}
\begin{equation}
\bo M'  + \bo d_3  \times \bo F= \rho ( I_2 \bo d_1 \times \ddot{\bo
d_1}+ I_1 \bo d_2 \times \ddot{\bo d_2}), \label{eq:amomentum2}
\end{equation}

\noindent where $\rho$ is the (volumetric) mass density, $A$ the
cross-sectional area, $I_1$ and $I_2$ the second moment of area of
the cross-section about $\bo d_1$ and $\bo d_2$ respectively, and
$\dot{(~)}$ denotes differentiation with respect to time.

For a closed system of equations these balance equations need to be
supplemented by constitutive relations that characterise the
material behaviour of the rod. We assume the rod to be made of
homogeneous isotropic linear viscoelastic material so that
stress-strain relations, based on a model by Valverde et
al.~\cite{valverde-cosserat}, are
\begin{equation}
\begin{split}
M_1 = EI_1 (\kappa_1+\gamma_v \dot \kappa_1),\\
M_2 = EI_2 (\kappa_2+\gamma_v \dot \kappa_2), \\
M_3 = GJ (\kappa_3+\gamma_v \dot \kappa_3), \\
\end{split}
\label{eq:constitutive1}
\end{equation}
where $\kappa_1$ and $\kappa_2$ are the curvatures about $\bo d_1$
and $\bo d_2$, respectively, while $\kappa_3$ is the twist about
$\bo d_3$. The constant $\gamma_v$ is the viscoelastic coefficient
of the material, $E$ is Young's modulus, $G$ is the shear modulus
and $J$ is the second moment of area of the section about $\bo d_3$.
We shall assume that the section is symmetric with respect to the
principal axes, in which case $J=I_1+I_2$.

The $\kappa_i$ are the components of the curvature vector
\begin{equation}
\m{\bd$\kappa$}=\sum_{i=1}^{3} \kappa_{i}\,\m{\bd$d$}_i,
\label{eq:twist1}
\end{equation}
which governs the evolution in space of the frame of directors as
one moves along the centreline:
\begin{equation}
\m{\bd$d$}'_i=\m{\bd$\kappa$}\times \m{\bd$d$}_i \quad\quad
(i=1,2,3). \label{eq:twist}
\end{equation}

The constitutive relations (\ref{eq:constitutive1}) can be used to
replace the $\kappa_i$ in (\ref{eq:twist}) by moments, after which
the equations (\ref{eqtangent1}), (\ref{eq:lmomentum2}),
(\ref{eq:amomentum2}) and (\ref{eq:twist}) form a system of 18
differential equations for the 18
unknowns $(\bo x,\bo F,\bo M,\bo d_1,\bo d_2,\bo d_3)$.\\

\noindent {\bf Remark:} We ignore in this study secondary electrodynamic
effects (such as an induced emf and hence additional current in the conductor)
as a result of the motion of the wire in the magnetic field \cite{jackson}.
Since we are considering a steadily rotating wire these effects would be null
on the configuration of the wire. However, the same would not be true for the
stability analysis, which considers arbitrary time-dependent perturbations.
We assume that these induction effects are negligible.

\subsection{Equations of motion in a uniformly rotating frame}
We shall also be interested in steadily rotating solutions and therefore
we transform the equilibrium equations (\ref{eq:lmomentum2}) and
(\ref{eq:amomentum2}) to a coordinate frame $\{\bo e_1,\bo e_2,\bo
e_3\}$ that rotates with constant angular velocity $\bo
\omega=\omega \bo k$ about the $\bo k$ axis (and the axis of the rod
in its trivial unstressed state). Noting that the derivative with
respect to time of an arbitrary vector $\bo V(s,t)$ is given by
\begin{equation}
\left . \frac {d \bo V (s,t)}{d t}\right|_i=\left . \frac {d \bo V
(s,t)}{d t}\right|_m + \bo \omega \times \bo V (s,t),
\label{eq:derfixed}
\end{equation}
where $\left . \frac {d}{d t}\right|_i$ indicates the derivative
with respect to time in the inertial frame and $\left . \frac {d}{d
t}\right|_m$ stands for the derivative with respect to time in the
moving frame, the equations (\ref{eq:lmomentum2}) and (\ref{eq:amomentum2})
expressed relative to $\{\bo e_1,\bo e_2,\bo e_3\}$ become
\begin{equation}
\bo F' + IB_0\bo d_3 \times \bo e_3 = \rho A (\ddot{\bo x} +2 \bo
\omega \times \dot{\bo x} + \bo \omega \times (\bo \omega \times \bo
x)),
\label{eq:lmomentum4}
\end{equation}
\begin{equation}
\begin{split}
\bo M'  + \bo d_3  \times \bo F &= \rho I_2 (\bo d_1 \times
\ddot{\bo d_1}+ 2\bo d_1 \times (\bo \omega \times \dot{\bo
d_1})+(\bo \omega \cdot
\bo d_1)(\bo d_1 \times \bo \omega)) \\
&+ \rho I_1 (\bo d_2 \times \ddot{\bo d_2} + 2\bo d_2 \times (\bo
\omega \times \dot{\bo d_2})+(\bo \omega \cdot \bo d_2)(\bo d_2
\times \bo \omega)).
\label{eq:amomentum4}
\end{split}
\end{equation}
The second term on the right-hand side of equation (\ref{eq:lmomentum4})
is the Coriolis force, while the third term is the centrifugal force, as
a result of the rotating coordinate system.

Steadily rotating (whirling) solutions satisfy the equations
(\ref{eq:lmomentum4}) and (\ref{eq:amomentum4}) with the dotted
variables set to zero:
\begin{equation}
\bo F' + IB_0\bo d_3 \times \bo e_3 = \rho A \bo \omega \times (\bo
\omega \times \bo x),
\label{eq:lmomentum4_whirl}
\end{equation}
\begin{equation}
\bo M'  + \bo d_3  \times \bo F = \rho I_2 (\bo \omega \cdot\bo d_1)
(\bo d_1 \times \bo \omega) + \rho I_1 (\bo \omega \cdot \bo
d_2)(\bo d_2\times \bo \omega).
\label{eq:amomentum4_whirl}
\end{equation}
The other equations (\ref{eqtangent1}) and (\ref{eq:twist}) do not
change their form, but all vectors are now to be considered as
expressed relative to the rotating frame $\{\bo e_1,\bo e_2,\bo
e_3\}$. Statical solutions are simply obtained by setting $\omega$
equal to zero.

For a well-posed problem the final 18 ODEs require 18 boundary
conditions to be specified, which we do next.

\subsection{Coat hanger boundary conditions}
\label{coat_hanger}

Helical solutions in rods are usually studied in infinitely long
rods, which avoids the need for imposing boundary conditions.
Indeed, it is not immediately clear how an exact helix can be
supported: the boundary conditions cannot be simply welded as no two
points on a helix have coaxial tangents, nor can they be simply
pinned because a helix has curvature and therefore carries a bending
moment. Here we formulate a set of boundary conditions that support
exact helical solutions. We call them coat hanger boundary
conditions, for obvious reasons.

\begin{figure}
\begin{center}
\includegraphics[width=0.5\linewidth]{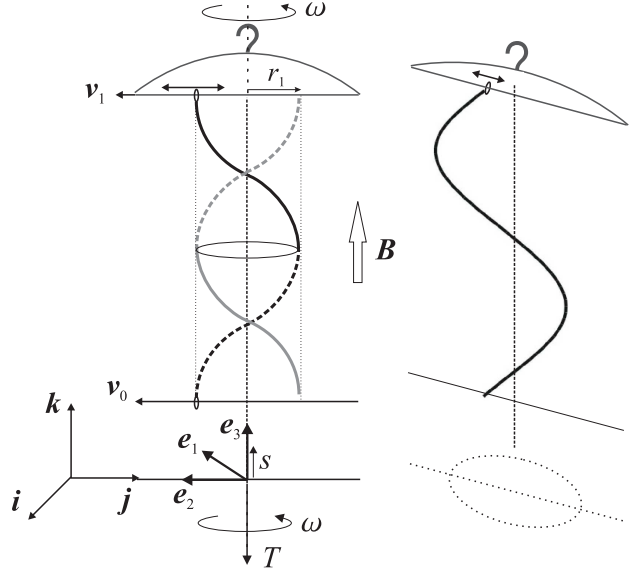}
\end{center}
\caption{Coat hanger boundary conditions.} \label{fig:hel:BC}
\end{figure}

Consider Fig.~\ref{fig:hel:BC} where a rod is suspended between
two axes $\bo v_0$ (at $s=0$) and $\bo v_1$ (at $s=1$) lying in two
parallel planes normal to $\bo e_3$. Axis $\bo v_1$ is taken to be
fixed in space, while $\bo v_0$ is free to move along $\bo e_3$. We
assume the axes to have a fixed relative rotation $\chi$, i.e., $\bo
v_0\cdot\bo v_1=\cos \chi$. The rod is free to hinge about and slide
along both $\bo v_0$ and $\bo v_1$. For definiteness we assume
that the rod is mounted in such a way that both axes $\bo v_0$ and
$\bo v_1$ are directed along the vector $\bo d_2$ in the rod's
cross-section. This situation is described by the following boundary
conditions:
\begin{eqnarray}
\bo d_1 (0,t) \cdot \bo v_0 &=& 0, \label{eq:BC36}\\
\bo d_3 (0,t) \cdot \bo v_0 &=& 0, \label{eq:BC37}\\
\bo M(0,t) \cdot \bo v_0 &=& 0, \label{eq:BC38}\\
\bo F(0,t) \cdot \bo v_0 &=& 0, \label{eq:BC39}\\
\bo F(0,t) \cdot \bo e_3 &=& T, \label{eq:BC40}\\
\bo x(0,t) \cdot (\bo v_0 \times \bo e_3) &=& 0, \label{eq:BC41}
\end{eqnarray}

\noindent at $s=0$, and
\begin{eqnarray}
\bo d_1 (L,t) \cdot \bo v_1 = 0, \label{eq:BC42}\\
\bo d_3 (L,t) \cdot \bo v_1 = 0, \label{eq:BC43}\\
\bo M(L,t) \cdot \bo v_1 = 0, \label{eq:BC44}\\
\bo F(L,t) \cdot \bo v_1 = 0, \label{eq:BC45}\\
\bo x(L,t) \cdot (\bo v_1 \times \bo e_3) = 0, \label{eq:BC46}\\
z(L,t) = L, \label{eq:BC47}
\end{eqnarray}
at $s=1$, where the position vector has been decomposed as
$\bo x = x \bo e_1 + y \bo e_2 + z \bo e_3$ and $T$ is an applied end
force (positive for tension). Conditions (\ref{eq:BC41}) and (\ref{eq:BC46})
restrict the movement of the ends of the rod to the planes spanned by
$(\bo v_0,\bo e_3)$ and $(\bo v_1,\bo e_3)$, respectively. To these 12
conditions we have to add conditions that ensure the orthonormality of the
director basis, for which we can take
\begin{equation}
\begin{split}
\bo d_1 (0,t) \cdot \bo d_1 (0,t) = 1, \\
\bo d_2 (0,t) \cdot \bo d_2 (0,t) = 1, \\
\bo d_3 (0,t) \cdot \bo d_3 (0,t) = 1, \\
\bo d_1 (0,t) \cdot \bo d_2 (0,t) = 0, \\
\bo d_1 (0,t) \cdot \bo d_3 (0,t) = 0, \\
\bo d_2 (0,t) \cdot \bo d_3 (0,t) = 0, \\
\end{split}
\label{eq:et:orthonorm}
\end{equation}
for a total of 18 boundary conditions, as required.

We shall take $\chi=0$ so that the initial rod, lying straight along
$\bo e_3$, is untwisted. This choice implies that any helical
solutions will have an integer number of half helical periods. Also,
the directors, and hence the cross-section of the rod, will make a
half-integer number of turns between $s=0$ and $s=L$. Without loss of
generality we may choose $\bo v_0=\bo v_1=\bo e_2$, so that at the
ends of the initial rod the directors $\{\bo d_1,\bo d_2,\bo d_3\}$ are
aligned with $\{\bo e_1,\bo e_2,\bo e_3\}$.

Since $\bo v_0$ and $\bo v_1$ are chosen aligned, there is a rigid-body
degree of freedom of translation of any solution along these axes. To
eliminate this degeneracy we replace condition (\ref{eq:BC45}) above by
\begin{equation}
\bo x(L,t) \cdot \bo v_1 = r_1, \label{eq:BC45a}
\end{equation}
where the slide $r_1$ along $\bo v_1$ (see Fig.~\ref{fig:hel:BC})
is chosen as follows. First note that a helix is a curve of
constant axial radius $r=\sqrt{x^2+y^2}$ and (total) curvature
$\kappa=\sqrt{\kappa_1^2+\kappa_2^2}$. The two are related by
$\kappa r=\sin^2\theta$, where $\theta$ is the helical angle
defined by $\bo d_3 \cdot \bo e_3=\cos\theta$ if $\bo e_3$ is along
the axis of the helix (the angle $\frac{\pi}{2}-\theta$ is usually
called the pitch angle). Since the rod is hinged about $\bo d_2$, we
have $\kappa_2(L,t)=0$. So, in order to ensure that any bifurcating
helix is centred at $\bo e_3$ we take
\begin{equation}
r_1=\frac{1-(\bo d_3(L,t) \cdot \bo e_3)^2}{\kappa_1(L,t)}, \label{eq:radius}
\end{equation}
giving a nonlinear boundary condition.
Note that $\kappa_1$ can here be taken with its sign, so that
(\ref{eq:BC45a}) also specifies which way the rod moves along $\bo
v_1$. When the rod buckles, the axis $\bo v_0$ lifts up and the rod
is free to find its own radius $r=|r_1|$. We stress that condition
(\ref{eq:BC45a}) has no effect on the bifurcation behaviour. In
particular, it does not suppress any non-helical solutions. It merely
ensures that if a helical solution bifurcates it will be centred at
the axis of rotation. This is important when we start rotating the axes
$\bo v_0$ and $\bo v_1$ about $\bo e_3$. A centred helix will experience
a uniform centrifugal force and is therefore expected to remain helical.

Of course the above coat hanger boundary conditions merely {\it
allow} for helical solutions. They need not exist. However, if the
equilibrium equations do have helical solutions and one-parameter
curves of such solutions intersect the trivial path of straight
solutions, then one might expect to detect them as (pitchfork)
bifurcations at critical buckling loads. The results presented in
Section \ref{results} show that this is indeed the case.

\section{Nondimensionalisation}

We make the system of equations dimensionless by scaling the
variables in the following way
\begin{equation}
\begin{split}
&\omega_c = \sqrt{\frac{EI_1}{\rho A L^4}},\quad \bar t=t\omega_c,
\quad \bar s=\frac{s}{L}\in[0,1], \quad \bar
\omega=\frac{\omega}{\omega_c}, \quad \bar{\bo x}=\frac{\bo x}{L},
\\
&\bar {\bo F}=\bo F \frac{L^2}{EI_1},\quad \bar{T}=T \frac{L^2}{E
I_1}, \quad \bar{\bo M}=\bo M\frac{L}{E I_1}, \quad \bar{\bo
\kappa}= \bo\kappa L.
\end{split}
\label{eq:nondim}
\end{equation}

\noindent Here $\omega_c$ is a reference characteristic bending
frequency of the rod.

With this nondimensionalisation the equations become (dropping the
overbars for simplicity and letting a prime denote $\frac{d}{d \bar s}$
and an overdot $\frac{d}{d \bar t}$):
\begin{eqnarray}
&\bo F' + B \bo d_3 \times \bo k =\ddot{\bo x} +2 \bo \omega \times
\dot{\bo x} + \bo \omega \times (\bo \omega
\times \bo x), \label{eq:lmomentum}\\
\bo M' + \bo d_3  \times \bo F &= P\left[R (\bo d_1 \times \ddot{\bo
d_1}+ 2\bo d_1 \times \bo \omega \times \dot{\bo d_1}+(\bo \omega
\cdot \bo d_1)(\bo d_1 \times \bo \omega)) \right. \nonumber \\
& \left. + (\bo d_2 \times \ddot{\bo d_2} + 2\bo d_2 \times \bo \omega \times
\dot{\bo d_2}+(\bo \omega \cdot \bo d_2)(\bo d_2 \times \bo
\omega)) \right], \label{eq:amomentum}\\
&\bo x'=\bo d_3, \label{eq:tangent}\\
&\bo d_i'=\bo \kappa \times \bo d_i, \label{eq:twist2}
\end{eqnarray}
\noindent and the constitutive relations can be written as
\begin{equation}
\bo M=(\kappa_1+\gamma \dot{\kappa_1})\bo d_1+
R(\kappa_2+\gamma \dot{\kappa_2})\bo d_2+ \frac{\Gamma
(1+R)}{2}(\kappa_3+\gamma \dot{\kappa_3})\bo d_3,
\label{eq:constitutive}
\end{equation}
\noindent where the dimensionless parameters are
\begin{equation}
P=\frac{I_1}{A L^2},\quad R=\frac{I_2}{I_1},\quad B=\frac{B_0 I
L^3}{E I_1},\quad \Gamma = \frac{2G}{E},\quad
\gamma=\gamma_v\omega_c ,
\end{equation}
and ($\frac{1}{\Gamma}-1$) is equal to Poisson's ratio. For the
boundary conditions we can still use (\ref{eq:BC36}) to
(\ref{eq:BC47}) if we assume that they now refer to
dimensionless variables and that the right-hand conditions are
imposed at $\bar{s}=1$.

\section{Perturbation scheme -- linearisation}

We consider whirling solutions (relative equilibria) that are
stationary in the moving frame $\{\bo e_1,\bo e_2,\bo e_3\}$. Such
solutions are found by solving the set of equations
(\ref{eq:lmomentum})--(\ref{eq:constitutive}) with the dotted
variables set to zero (thus obtaining an ODE). To study their
stability we linearise the full PDE
(\ref{eq:lmomentum})--(\ref{eq:constitutive}) about these whirling
solutions. We follow the approach in \cite{valverde-cosserat}, which
is similar to approaches in \cite{fraser,goriely}. The stability of
static (non-whirling) solutions can be investigated by simply
setting the angular velocity $\omega$ to zero.

We start our perturbation analysis by writing
\begin{equation}
\bo d_i(s,t)=\bo d_i^{0}(s)+\delta \bo
d_i^{t}(s,t)+O(\delta^2),\quad i=1,2,3, \label{per:dir}
\end{equation}
where $\bo d_i^{0}(s)$ are the basis vectors of a quasi-stationary solution,
$\bo d_i^{t}(s,t)$ are the basis vectors of a time-dependent perturbation
and $\delta$ is a small bookkeeping parameter introduced to separate scales.
Note that, in order to preserve orthonormality to $O(\delta)$ ($\bo d_i
\cdot \bo d_j = \delta_{ij}+O(\delta^2)$), we must have
\begin{equation}
\bo d_i^{t}(s,t)=\sum_{j=1}^3 A_{ij}(s,t) \bo d_j^0(s),\quad
i=1,2,3, \label{per:dirt}
\end{equation}
where the matrix $A_{ij}$ is skew-symmetric and can be written as
\begin{equation}
%\begin{displaymath}
\bo A= \left( \begin{array}{ccc}
0 & \alpha_3 & -\alpha_2 \\
-\alpha_3 & 0 & \alpha_1 \\
\alpha_2 & -\alpha_1  & 0
\end{array} \right).
\label{per:dirt2}
%\end{displaymath}
\end{equation}
Thus, the nine components of the director basis perturbation are
described by only three independent parameters, and if we introduce
\begin{equation}
\bo \alpha = (\alpha_1,\alpha_2,\alpha_3)^T
\label{alpha}
\end{equation}
(with respect to the unperturbed director basis) then the perturbed director
basis can be expressed as
\begin{equation}
\bo d_i(s,t)=\bo d_i^{0}(s)+\delta \bo \alpha(s,t) \times \bo
d_i^{0}(s) +O(\delta^2), \quad i=1,2,3. \label{per:dirfinal}
\end{equation}

Using (\ref{per:dirfinal}), the perturbation of an arbitrary vector
$\bo V=\sum_{i=1}^3 V_i \bo d_i$ can be written on the basis
$\{\bo d_1^0,\bo d_2^0,\bo d_3^0\}$ as
\begin{equation}
\bo V=\bo V^0+\delta \bo V^t  +O(\delta^2)
=\sum_{i=1}^3 [V_i^0+\delta( V_i^t +(\bo \alpha \times \bo
V^0)_i)] \bo d_i^0    +O(\delta^2), \label{per:varfinal}
\end{equation}
where $()_i$ denotes the component along $\bo d_i^0$ and time and space
dependence of the variables have been suppressed for the sake of simplicity
\cite{valverde-cosserat}.

Applying this perturbation scheme to the PDEs
(\ref{eq:lmomentum})--(\ref{eq:constitutive}) and the boundary conditions,
we arrive at an $O(1)$ nonlinear ODE for the quasi-stationary solutions and
an $O(\delta)$ linear PDE governing their stability.

\subsection{The $O(1)$ equations -- quasi-stationary whirl}
\label{stationaryBVP}

The $O(1)$ equations are time-independent. Recalling that $\bo
\omega= \omega \bo e_3$, we find the $O(1)$ terms of the linear
momentum equation (\ref{eq:lmomentum}), projected on the director
basis $\{\bo d_1^0,\bo d_2^0,\bo d_3^0\}$, to give

{\setlength\arraycolsep{2pt}
\begin{eqnarray}
(F_1^0)'-F_2^0\kappa_3^0+F_3^0\kappa_2^0 + B(d_{32}^0d_{11}^0-d_{31}^0
d_{12}^0)&=&-\omega^2(x^0d_{11}^0+y^0d_{12}^0),\label{stat:forces1}\\
(F_2^0)'-F_3^0\kappa_1^0+F_1^0\kappa_3^0 + B(d_{32}^0d_{21}^0-d_{31}^0
d_{22}^0)&=&-\omega^2(x^0d_{21}^0+y^0d_{22}^0),\label{stat:forces2}\\
(F_3^0)'-F_1^0\kappa_2^0+F_2^0\kappa_1^0&=&-\omega^2(x^0d_{31}^0+y^0d_{32}^0),
\label{stat:forces3}
\end{eqnarray}}

\noindent where subscripts are used to indicate components relative
to the basis vectors $\{\bo d_1^0,\bo d_2^0,\bo d_3^0\}$ (but the
$\bo d_i^0$ components are relative to $\{\bo e_1,\bo e_2,\bo
e_3\}$). Similarly, the $O(1)$ term of the angular momentum equation
(\ref{eq:amomentum}), projected on the director basis $\{\bo
d_1^0,\bo d_2^0,\bo d_3^0\}$ gives

{\setlength\arraycolsep{2pt}
\begin{eqnarray}
(M_1^0)'&=&\frac{2
M_3^0M_2^0}{\Gamma(1+R)}-\frac{M_2^0M_3^0}{R}+F_2^0
+P\omega^2d_{23}^0(d_{22}^0d_{11}^0-d_{21}^0d_{12}^0),\label{stat:moments1}\\
(M_2^0)'&=&-\frac{2 M_3^0M_1^0}{\Gamma(1+R)}+M_1^0M_3^0-F_1^0
+PR\omega^2d_{13}^0(d_{21}^0d_{12}^0-d_{11}^0d_{22}^0),\label{stat:moments2}\\
(M_3^0)'&=&\frac{M_2^0M_1^0}{R}-M_1^0M_2^0+PR\omega^2d_{13}^0
(d_{12}^0d_{31}^0-d_{11}^0d_{32}^0)+P\omega^2d_{23}^0(d_{22}^0d_{31}^0-d_{21}^0d_{32}^0).\label{stat:moments3}
\end{eqnarray}}

The $O(1)$ term of equation (\ref{eq:tangent}) can be expressed as
\begin{equation}
(\bo x^0)'=\bo d_{3}^0, \label{stat:tangent}
\end{equation}
and the twist equation (\ref{eq:twist2}) by
\begin{equation}
 (\bo d_{i}^0)'=\bo \kappa^0 \times \bo d_{i}^0, \qquad i=1,2,3,
\label{stat:twist}
\end{equation}
where $\bo \kappa^0 = \sum_{j=1}^3 \kappa_j^0 \bo d_{j}^0$. The
$O(1)$ term of constitutive relations (\ref{eq:constitutive}) can be
expressed as
\begin{equation}
M_1^0=\kappa_1^0, \quad\quad M_2^0=R\kappa_2^0, \quad\quad
M_3^0=\frac{\Gamma(1+R)}{2}\kappa_3^0, \label{stat:constitutive}
\end{equation}
which can be used to express the $\kappa_i^0$ in (\ref{stat:twist})
in terms of the moments $M_i^0$.

Proceeding in the same way, the $O(1)$ part of the boundary
conditions is given by
\begin{eqnarray}
\bo d_1^0 (0) \cdot \bo v_0 &=& 0, \label{eq:stat:BC36}\\
\bo d_3^0 (0) \cdot \bo v_0 &=& 0, \label{eq:stat:BC37}\\
\bo M^0(0) \cdot \bo v_0 &=& 0, \label{eq:stat:BC38}\\
\bo F^0(0) \cdot \bo v_0 &=& 0, \label{eq:stat:BC39}\\
\bo F^0(0) \cdot \bo e_3 &=& T, \label{eq:stat:BC40}\\
\bo x^0(0) \cdot (\bo v_0 \times \bo e_3) &=& 0,
\label{eq:stat:BC41}
\end{eqnarray}
\begin{eqnarray}
\bo d_1^0 (1) \cdot \bo v_1 &=& 0, \label{eq:stat:BC42}\\
\bo d_3^0 (1) \cdot \bo v_1 &=& 0, \label{eq:stat:BC43}\\
\bo M^0(1) \cdot \bo v_1 &=& 0, \label{eq:stat:BC44}\\
\bo x^0(1) \cdot (\bo v_1 \times \bo e_3) &=& 0, \label{eq:stat:BC45}\\
\bo x^0(1) \cdot \bo v_1 &=& r^0(1) = \frac{1}{\kappa_1^0(1)}\left(
1-(d_{33}^0(1))^2\right), \label{eq:stat:BC46}\\
z^0(1) &=& 1, \label{eq:stat:BC47}
\end{eqnarray}
where (\ref{eq:stat:BC46}) is the $O(1)$ contribution from (\ref{eq:BC45a}),
with the radius of the helix also affected by the perturbation scheme, i.e.,
$r(s,t)=r^0(s)+\delta r^t(s,t)$.

\subsection{The $O(\delta)$ equations -- linearisation}

The $O(\delta)$ part of the linear momentum equation (\ref{eq:lmomentum})
can be written as
\begin{equation}
\begin{split}
(\bo F^t(s,t))'+& \bo B_1(s)\bo F^t(s,t)+\bo B_2(s)\bo x^t(s,t)+\bo
B_3(s)\bo \alpha'(s,t)+\bo B_4(s)\bo \alpha(s,t)\\
=&\bo B_5(s) \ddot{\bo x}^t(s,t)+\bo B_6(s) \dot{\bo x}^t(s,t),
\end{split} \label{per:lmomentum}
\end{equation}
where the $3\times 3$ matrices $\bo B_i(s)$ are given in
Appendix A. Here we have expressed $\bo F^t$ relative to
$\{\bo d_1^0,\bo d_2^0,\bo d_3^0\}$ and $\bo x^t$ relative to
$\{\bo e_1,\bo e_2,\bo e_3\}$. For the $O(\delta)$ part of the
angular momentum equation (\ref{eq:amomentum}) we can write
\begin{equation}
\begin{split}
(\bo M^t(s,t))'+& \bo C_1(s)\bo M^t(s,t)+\bo C_2(s)\bo
\alpha'(s,t)+\bo C_3(s)\bo \alpha(s,t)+ \bo C_4(s)\bo F^t(s,t)\\
=& \bo C_5(s) \ddot{\bo \alpha}(s,t)+\bo C_6(s) \dot{\bo
\alpha}(s,t),
\end{split} \label{per:amomentum}
\end{equation}
where the matrices $\bo C_i(s)$ are again given in Appendix A.
$\bo M^t$ is expressed relative to $\{\bo d_1^0,\bo d_2^0,\bo d_3^0\}$.
The 9 twist equations (\ref{eq:twist2}) at $O(\delta)$ are reduced to
only 3  independent equations that relate $\bo \kappa^t$ and $\bo
\alpha$ as
\begin{equation}
\bo \kappa^t(s,t)=\bo \alpha'(s,t) +\bo \kappa^0(s) \times \bo
\alpha(s,t). \label{per:twist}
\end{equation}
Introducing these relations into the $O(\delta)$ part of the
constitutive relations gives
\begin{equation}
\bo M^t(s,t) + \bo D_1(s)\bo \alpha'(s,t)+\bo D_2(s)\bo \alpha(s,t)
= \bo D_3(s) \dot{\bo \alpha}(s,t)+\bo D_4(s)
\frac{\partial^2}{\partial s \partial t}(\bo \alpha(s,t)),
\label{per:constitutive}
\end{equation}
where the matrices $\bo D_i(s)$ are given in Appendix A. Finally, the
$O(\delta)$ part of equation (\ref{eq:tangent}) yields
\begin{equation}
\begin{split}
(\bo x^t(s,t))' = \alpha(s,t)\times\bo d_3^0(s).
\label{per:tangent}
\end{split}
\end{equation}

Applying the perturbation scheme to the boundary conditions at $O(\delta)$,
we obtain
{\setlength\arraycolsep{2pt}
\begin{eqnarray}
\alpha_1(0,t)&=& 0,\label{per:bc1}\\
\alpha_3(0,t)&=& 0,\label{per:bc2}\\
x^t(0,t)&=&0,\label{per:bc3}\\
d_{13}^0(0)F_1^t(0,t)+d_{33}^0(0)F_3^t(0,t)+(d_{13}^0(0)F_3^0(0)-d_{33}^0(0)F_1^0(0))\alpha_2(0,t)&=&0, \label{eq:bc4}\\
F_2^t(0,t)&=& 0,\label{per:bc5}\\
M_2^t(0,t)&=& 0,\label{per:bc6}\\
x^t(1,t)&=& 0,\label{per:bc7}\\
z^t(1,t)&=& 0,\label{per:bc8}\\
\alpha_1(1,t)&=& 0,\label{per:bc9}\\
\alpha_3(1,t)&=& 0,\label{per:bc10}\\
M_2^t(1,t)&=& 0,\label{per:bc11}\\
y^t(1,t)+\frac{2}{\kappa_1^0(1)}(\bo d_3^0(1)\cdot\bo e_3)
(\bo d_3^t(1,t)\cdot\bo e_3) + \frac{1}{{\kappa_1^0(1)}^2}\left(1-
(\bo d_3^0(1)\cdot\bo e_3)^2\right)\kappa_1^t(1,t)&=&0.\label{per:bc12}
\end{eqnarray}}

After elimination of the $\kappa_i^0$ by means of (\ref{stat:constitutive}),
the set of 12 equations (\ref{per:lmomentum}), (\ref{per:amomentum}),
(\ref{per:constitutive}) and (\ref{per:tangent}) together with the 12
boundary conditions (\ref{per:bc1})--(\ref{per:bc12}), with appropriate
initial conditions form a well-posed initial-boundary-value problem.

\section{Stability analysis}
\label{sect:stab_anal}

Since we are interested in stability of solutions we look for
solutions of the $O(\delta)$ equations of the form
\begin{eqnarray}
\bo x^t (s,t) = \hat{\bo x}^t (s) e^{\lambda t},\label{per:sol1}\\
\bo \alpha (s,t) = \hat{\bo \alpha} (s) e^{\lambda t},\\
\bo F^t (s,t) = \hat{\bo F}^t (s) e^{\lambda t},\\
\bo M^t (s,t) = \hat{\bo M}^t (s) e^{\lambda t}.\label{per:sol4}
\end{eqnarray}
When these expressions are inserted into
(\ref{per:lmomentum})--(\ref{per:tangent}) a linear eigenvalue
problem for a 12-dimensional ODE is obtained in terms of the
variables $(\hat{\bo x}^t,\hat{\bo \alpha}^t,\hat{\bo F}^t,\hat{\bo
M}^t)$. The eigenvalue $\lambda$ measures the growth of small
perturbations and is to be found as part of the solution.
Eigenvalues come as complex conjugate pairs. A whirling state is
unstable if at least one of the (in general infinitely many)
$\lambda$'s has positive real part.

To solve a real system of equations we split the eigenvalues and variables
(eigenfunctions) into real and imaginary parts,
$\lambda=\lambda_r + i \lambda_i$, $ \hat{\bo x}^t = \hat{\bo x}^t_r
+ i \hat{\bo x}^t_i$, $ \hat{\bo \alpha}^t = \hat{\bo \alpha}^t_r +
i \hat{\bo \alpha}^t_i$, $ \hat{\bo F}^t = \hat{\bo F}^t_r + i
\hat{\bo F}^t_i$ and $ \hat{\bo M}^t = \hat{\bo M}^t_r + i \hat{\bo
M}^t_i$. The equations (\ref{per:lmomentum})--(\ref{per:tangent}) along with
the boundary conditions (\ref{per:bc1})--(\ref{per:bc12}), are similarly
split into real and imaginary parts. Thus we end up with a doubled
24-dimensional linearised boundary-value problem.

\subsection{Stability of the straight rod -- static magnetic buckling}
\label{sect:magn_buckling}

The trivial solution of the $O(1)$ equations
(\ref{stat:forces1})--(\ref{stat:constitutive}), representing a straight and
untwisted rod, is given by
\begin{equation}
\bo x(s)=s\bo e_3, \quad \bo F(s)=-T\bo e_3, \quad \bo M(s)=\bo
0, \quad \bo d_i(s)=\bo e_i~~~(i=1,2,3), \quad s\in [0,1].
\label{trivial_sol}
\end{equation}
It satisfies the coat hanger boundary conditions.
For the statics case ($\omega=0$) without end force ($T=0$) the $O(\delta)$
equations (\ref{per:lmomentum}), (\ref{per:amomentum}),
(\ref{per:constitutive}), (\ref{per:tangent}) about this trivial solution,
on inserting (\ref{per:sol1})--(\ref{per:sol4}), can be written as
\begin{eqnarray}
&& x''''-\lambda^2Px''+\frac{\lambda^2}{R}x-\frac{B}{R}y'=0, \nonumber \\
&& y''''-\lambda^2Py''+\lambda^2y+Bx'=0, \label{lin_eqs} \\
&& M_3''-\frac{2\lambda^2P}{\Gamma}M_3=0, \nonumber
\end{eqnarray}
with boundary conditions
\begin{equation}
x(0)=x(1)=x''(0)=x''(1)=y'''(0)=y(1)=y'(0)=y'(1)=0, \quad\quad
M_3'(0)=0=M_3'(1),
\label{BCs}
\end{equation}
while $F_3\equiv 0$, $z\equiv 0$. Note that the torsional ($M_3$) modes
decouple from the bending ($x,y$) modes.

To find the static magnetic buckling loads we set $\lambda=0$. The bending
equations then reduce to
\begin{equation}
z''''''+\frac{B^2}{R}z=0, \quad\quad \mbox{for} \quad\quad z=x',
\end{equation}
subject to
\begin{equation}
x(0)=x(1)=x''(0)=x''(1)=x''''(0)=x''''(1)=0.
\label{BC_x}
\end{equation}
On setting $z=e^{iks}$ we obtain the characteristic equation $-k^6+B^2/R=0$
with solutions
$k_{1,2}=\pm\beta,~k_{3,4,5,6}=\pm\beta\left(1\pm i\sqrt{3}\right)/2$,
where $\beta=B^{1/3}/R^{1/6}$. Application of the boundary conditions
(\ref{BC_x}) to the general solution $z(s)=\sum_{j=1}^6 a_je^{ik_js}$ leads
to the remarkably simple condition:
\begin{equation}
B=(n\pi)^3\sqrt{R}.
\label{char_coat}
\end{equation}
These critical loads correspond to pitchfork bifurcations where non-trivial
solutions bifurcate from the trivial straight solution. We stress that the
above calculation is only possible for the statics case. If $\omega\neq 0$
then the $x$ and $y$ equations do not decouple and no simple characteristic
equation is obtained. However, bifurcating branches of helical solutions,
and hence critical loads, can be computed explicitly, even for non-zero
$\omega$; see Appendix B.

\subsection{Eigenvalues for the unperturbed problem ($T=0$, $\gamma=0$,
$\omega=0$, $B=0$)}
\label{sect:unperturbed}

We shall call the case where $T=0$, $\gamma=0$, $\omega=0$ and $B=0$ the
{\it unperturbed problem}. For this problem explicit expressions can be
obtained for the eigenvalues of the linearisation about the straight solution.
The $x$ and $y$ equations in (\ref{lin_eqs}) decouple into two fourth-order
beam equations:
\begin{equation}
\begin{split}
x''''-\lambda^2Px''+\frac{\lambda^2}{R}x=0, \vspace*{0.2cm} \\
y''''-\lambda^2Py''+\lambda^2y=0,
\end{split}
\label{lin_eqs_0}
\end{equation}
subject to boundary conditions (\ref{BCs}). Since we anticipate imaginary
eigenvalues we set $\lambda=i\mu$, $x=e^{iks}$, $y=e^{i\kappa s}$, and find
for the $x$ equation
\begin{eqnarray}
&& k_{1,2}=\pm\left(\frac{1}{2}\mu^2P+\frac{1}{2}
\sqrt{\mu^4P^2+4\mu^2/R}\right)^{1/2}=:\pm a, \nonumber \\
&& k_{3,4}=\pm i\left(\frac{1}{2}\sqrt{\mu^4P^2+4\mu^2/R}-
\frac{1}{2}\mu^2P\right)^{1/2}=:\pm i b, \nonumber
\end{eqnarray}
while for the $y$ equation
\begin{eqnarray}
&& \kappa_{1,2}=\pm\left(\frac{1}{2}\mu^2P+\frac{1}{2}
\sqrt{\mu^4P^2+4\mu^2}\right)^{1/2}=:\pm\alpha, \nonumber \\
&& \kappa_{3,4}=\pm i\left(\frac{1}{2}\sqrt{\mu^4P^2+4\mu^2}-
\frac{1}{2}\mu^2P\right)^{1/2}=:\pm i\beta, \nonumber
\end{eqnarray}
where $a$, $b$, $\alpha$, $\beta$ are non-negative real numbers. The general
solutions are
\begin{eqnarray}
&& x(s)=A_x\sin as + B_x\cos as + C_x\sinh bs + D_x \cosh bs, \nonumber \\
&& y(s)=A_y\sin \alpha s + B_y\cos \alpha s + C_y\sinh \beta s +
D_y \cosh \beta s. \nonumber
\end{eqnarray}
Application of the boundary conditions (\ref{BCs}) leads to
\[
\sin a=0 \quad\quad \mbox{and} \quad\quad
\beta\cos\alpha\sinh\beta + \alpha\sin\alpha\cosh\beta=0.
\]
The first equation implies
\[
\mu=\pm\frac{n^2\pi^2}{\sqrt{n^2\pi^2P+1/R}}, \quad n=1,2,3,....
\]
The second equation is transcendental and needs to be solved numerically
to obtain the eigenvalues $\mu$, for instance by using a Newton-Raphson
scheme. Meanwhile, the torsional eigenvalues for the $M_3$ equation in
(\ref{lin_eqs}) are given by
\[
\mu=\pm n\pi\sqrt{\frac{\Gamma}{2P}}, \quad n=1,2,3,....
\]
These are all the eigenvalues for the unperturbed problem. They will be used
as starting values in the numerical procedure described next.

\subsection{Numerical procedure}

The main idea is to use the known eigenvalues in the unperturbed problem as
starting values in a continuation procedure in order to compute the
eigenvalues and corresponding eigenfunctions for general values of the
parameters $T$, $\gamma$, $\omega$ and $B$. For this we use the well-tested
code AUTO \cite{doedel} (specifically AUTO2000). AUTO solves boundary-value
problems by means of orthogonal collocation. It requires a starting solution
and can then trace out solution curves as a parameter of the problem is
varied. Bifurcations are detected where branches of solutions intersect.
At such points AUTO is able to switch branches and compute curves of
bifurcating solutions.

Our procedure takes advantage of the fact that $\lambda$ appears only
quadratically in the linearisation (\ref{per:lmomentum}),
(\ref{per:amomentum}), (\ref{per:constitutive}) and (\ref{per:tangent}) if
$\gamma=0$ and $\omega=0$. To explain the method consider the typical
$O(\delta)$ equation
\begin{equation}
z''''-\lambda^2f(s)z''+\lambda^2g(s)z=0,
\end{equation}
where $f$ and $g$ are functions of the $O(1)$ solution. Writing $z=x+iy$,
$\lambda=\lambda_r+i\lambda_i$, we can decompose the $z$ equation into
\begin{equation}
\begin{split}
x''''-(\lambda_r^2-\lambda_i^2)f(s)x''+2\lambda_r\lambda_i f(s)y''+
(\lambda_r^2-\lambda_i^2)g(s)x-2\lambda_r\lambda_i g(s)y=0, \\
y''''-(\lambda_r^2-\lambda_i^2)f(s)y''-2\lambda_r\lambda_i f(s)x''+
(\lambda_r^2-\lambda_i^2)g(s)y+2\lambda_r\lambda_i g(s)x=0.
\end{split}
\label{re_im_eq}
\end{equation}
The important thing to note here is that these equations decouple
into two identical equations if the eigenvalue is either imaginary
($\lambda_r=0$) or real ($\lambda_i=0$).

This suggests the following sequence of steps, involving boundary-value
problem of increasing dimension, to compute eigenvalues of static or
uniformly whirling solutions.

\begin{enumerate}

\item
Consider the unperturbed problem of Section~\ref{sect:unperturbed} and,
noting that all eigenvalues are purely imaginary, solve the 30-dimensional
system of 18 $O(1)$ equations and one 12-dimensional system for the
imaginary part of the $O(\delta)$ equations (cf. the $y$ equation in
(\ref{re_im_eq})). Set $\lambda_r=0$ and use $\lambda_i$ as the continuation
parameter in AUTO in order to compute the eigenvalues (instead of solving the
transcendental equations in Section~\ref{sect:unperturbed}). These eigenvalues
will show up as branching points (BP), or pitchfork bifurcations, as
eigenvalues by definition are those values for which non-zero BVP solutions
exist. By symmetry it is only necessary to consider $\lambda_i>0$.

\item
Keeping the same 30-dimensional system, switch branches at a BP to
compute (`grow') the corresponding (imaginary) eigenfunction. Since the
equations are linear the value of $\lambda_i$ will not change in this run.
For later use we monitor the non-zero solution by means of some measure
$||.||_i$ (not necesarily a proper norm) on the space
$\{\hat{\bo x}_i^t,\hat{\bo \alpha}_i^t,\hat{\bo F}_i^t,\hat{\bo M}_i^t\}$
of imaginary linearised variables.

\item
Now consider the full system of 42 equations (18 $O(1)$ equations and two
sets of 12-dimensional $O(\delta)$ equations (cf.~(\ref{re_im_eq})). Fix the
measure $||.||_i$ on the imaginary part and release $\lambda_r$ instead in
order to compute the real eigenfunction (since the imaginary part of the
solution is fixed there is only one branch of solutions through the starting
point and there is nowhere else to go for the continuation but to `grow' the
real eigenfunction). Again we monitor this function by means of a suitable
measure $||.||_r$. In this run neither $\lambda_r$ nor $\lambda_i$ will
change.

This approach works because the solution obtained in step 2 also solves
the full 42-dimensional system when the extra 12 variables
$(\hat{\bo x}_r^t,\hat{\bo \alpha}_r^t,\hat{\bo F}_r^t,\hat{\bo M}_r^t)$
are set to zero. This is a consequence of the fact that the real and imaginary
parts of the $O(\delta)$ equations decouple if $\lambda_r=0$, as a result
of the quadratic dependence of the eigenvalue problem on $\lambda$
(cf.~(\ref{re_im_eq})).

\end{enumerate}

Steps 2 and 3 can be performed for as many of the BPs computed in
step 1 as required and will give the corresponding eigenvalues and
eigenfunctions. Once these have been obtained both measures
$||.||_r$ and $||.||_i$ can be fixed and an extra system parameter
such as $B$ or $\omega$ released in order to trace the eigenvalues
(and hence monitor stability changes) as system parameters are
varied. (Note that fixing $||.||_r$ and $||.||_i$ makes sense as
eigenfunctions are only defined up to a multiplicative factor.)

The above 3-step procedure is not limited to linearisations about the trivial
straight solution. It can be applied to any {\it starting solution} that has
no eigenvalue with both $\lambda_r$ and $\lambda_i$ non-zero, as these would
not be picked up in step 1. (It is of course no problem if eigenvalues
become fully complex (Hopf bifurcation) in the course of further
continuations.) For instance, we find that at the first critical $B$, given
by (\ref{char_coat}), the lowest conjugate pair of eigenvalues $\pm\lambda_i$
goes to zero and becomes a real pair of eigenvalues, signalling a stability
change of the straight rod. The (first-mode) solution bifurcating at this
point is stable with all eigenvalues being imaginary and the above procedure
can be applied to find the eigenvalues.

We end this section with a few comments:

\begin{enumerate}

\item[$(i)$]
There are infinitely many eigenvalues and the above procedure only finds the
lowest order ones. This is of course a limitation of any numerical scheme.
We find that eigenvalues vary slowly with system parameters, suggesting
that stability is governed by the lowest-order eigenvalues. We typically
consider 5 or 6 eigenvalues.

\item[$(ii)$]
Note that in steps 1 and 2 above we could not have taken the full
42-dimensional system of equations as that would have made the branching
points (pitchfork bifurcations) degenerate and AUTO would not detect a BP.
This is because if $\lambda_r=0$ (or $\lambda_i=0$) the two sets of
12-dimensional linearised equations are identical (cf.~(\ref{re_im_eq})).

\end{enumerate}

\section{Numerical results}\label{results}

\subsection{The statics case ($\omega=0$)} \label{coat_statics}

\begin{figure}
\begin{center}
\includegraphics[width=0.6\linewidth]{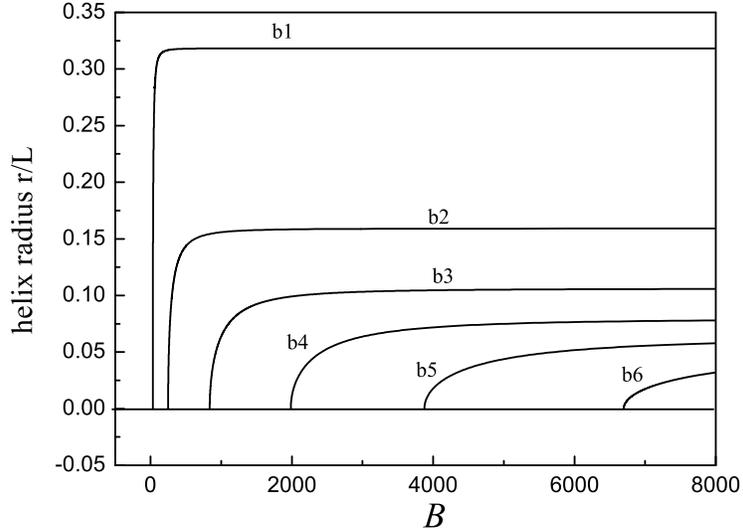}
\end{center}
\caption{Bifurcation diagram for a stationary isotropic rod subject to coat
hanger boundary conditions ($\omega=0$).}
\label{fig:hel:bif-helix}
\end{figure}

Fig.~\ref{fig:hel:bif-helix} shows the bifurcation diagram obtained
by varying the magnetic field parameter $B$ in the isotropic case
($R=1$). Throughout this entire section the dimensionless parameters
taken, unless stated otherwise, are those listed in
Table~\ref{table:hel:parameters}, where realistic dimensional
parameters are also given. Pitchfork bifurcations on the trivial branch
occur at $BP1=31.01$, $BP2=248.05$, $BP3=837.17$, etc., in agreement
with (\ref{char_coat}) (only the positive-$r$ branch is shown). The
bifurcating solutions are found to be exact helices and therefore the
dimensionless helical radius $r/L$ is used as solution measure on the
vertical axis. A slight complication in computing this diagram occurs
because of the denominator in (\ref{eq:radius}), which is zero for the
straight rod. However, this problem is easily resolved by replacing
boundary condition (\ref{eq:BC45a}) by $\bo x(L)\cdot\bo v_1=0$ along
the trivial branch and switching back to (\ref{eq:BC45a}) once an
incipient non-trivial solution has been obtained.

\begin{table}
  \centering
  \caption{Parameters used for the coat hanger boundary conditions.}
  \vspace{0.05\linewidth}
  \begin{tabular}{|c|c||c|c|} \hline
    $L$ & 5 m & $P$ & 0.001\\ \hline
    $A$ & $3 \times 10^{-5}$ m$^2$ & $R$ & 1\\ \hline
    $E$ & $30 \times 10^{9}$ N/m$^2$ & $\Gamma$ & 0.76923\\ \hline
    $I_1 = I_2$ & $9 \times 10^{-11}$ m$^4$ &  & \\ \hline
  \end{tabular}
  \label{table:hel:parameters}
\end{table}
\begin{figure}
\begin{center}
\includegraphics[width=0.45\linewidth]{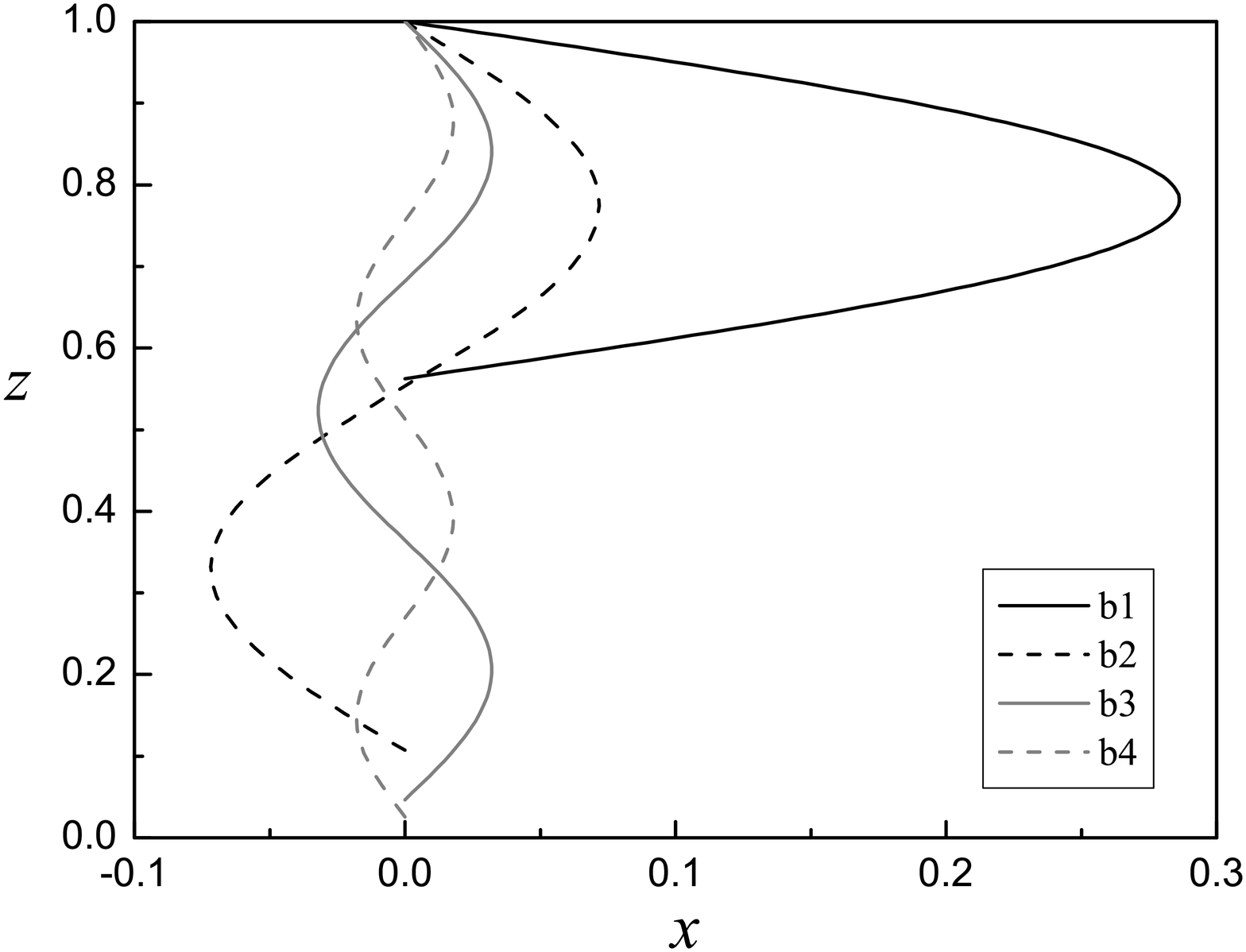}
~~~~~\includegraphics[width=0.45\linewidth]{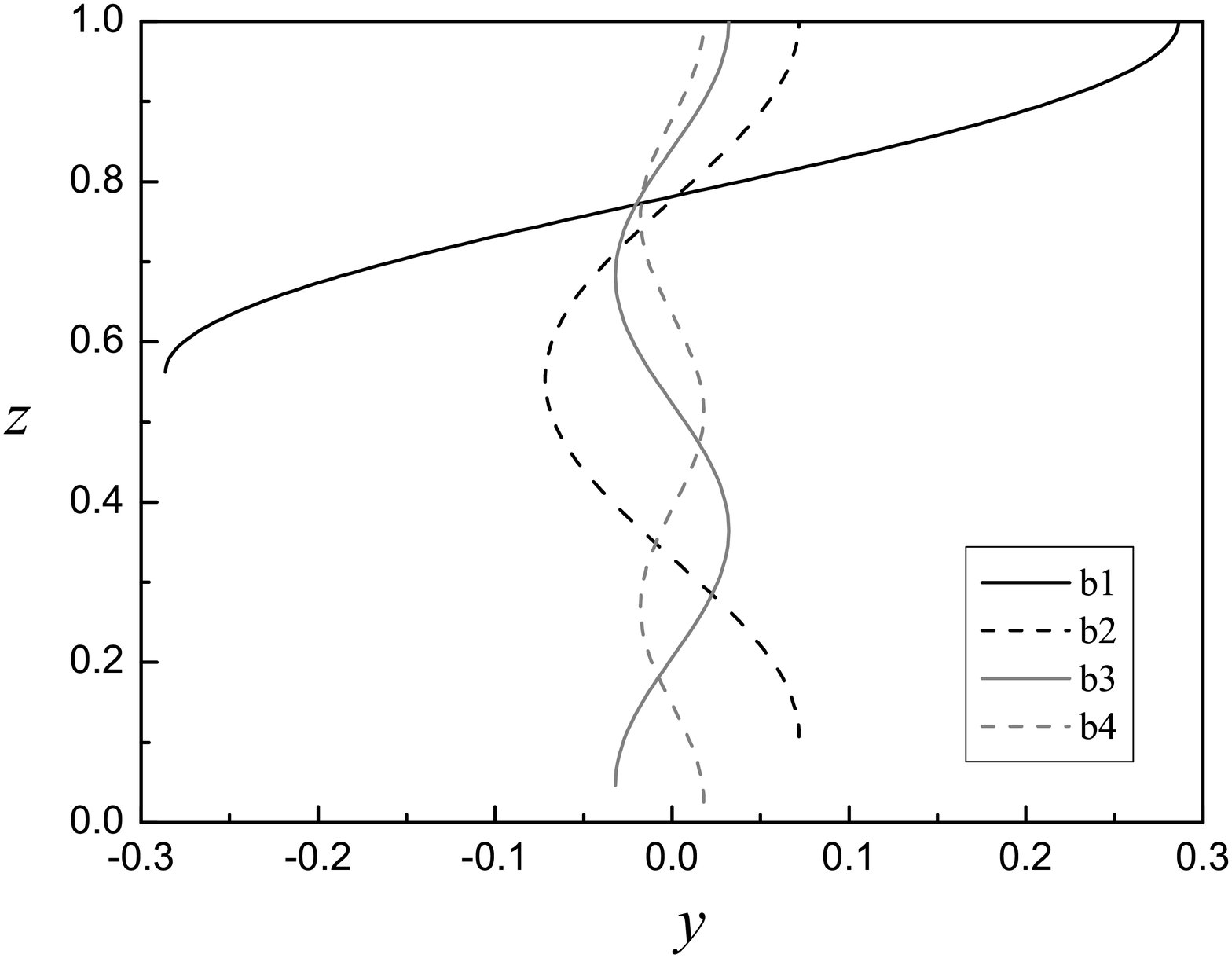}
\end{center}
\caption{The \{$\bo e_1$-$\bo e_3$\} and \{$\bo e_2$-$\bo e_3$\}
projections of the first 4 helical modes at constant curvature
$\kappa = 8$. Values of $B$ are: 61.68 (b1), 270.69 (b2), 868.72
(b3) and 2024.70 (b4).}
\label{fig:hel:hel-modes-shape}
\end{figure}
\begin{figure}
\begin{center}
\includegraphics[width=0.35\linewidth]{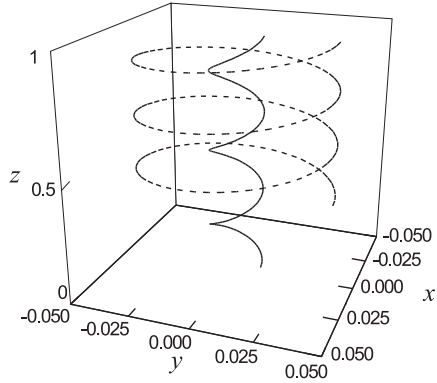}
\end{center}
\caption{3D view of two helices on branch b6.}
\label{fig:hel:hel-shape-3D}
\end{figure}

Fig.~\ref{fig:hel:hel-modes-shape} shows \{$\bo e_1$-$\bo e_3$\}
and \{$\bo e_2$-$\bo e_3$\} projections of bifurcating solutions
along the first four branches, taken at constant curvature
$\kappa=8$. It was noted in Section~\ref{coat_hanger} that the coat
hanger boundary conditions with $\chi=0$ (i.e., parallel end
supports $\bo v_0$ and $\bo v_1$) only allow helices of a
half-integer number, $n$, of helical turns. We find that each
successive bifurcating solution in Fig.~\ref{fig:hel:bif-helix}
has one more half helical turn. For large $B$ the solutions
approach a circular shape in the $z=1$ plane, with corresponding
value $r/L=1/(n\pi)$ along the vertical axis in
Fig.~\ref{fig:hel:bif-helix}. The bifurcating branches have handedness.
That is, regardless of the sign of $r$, the bifurcating solutions are
right-handed helices if, as here, $B>0$ and would be left-handed
helices if $B<0$, i.e., had we run $B$ in the other direction.
Fig.~\ref{fig:hel:hel-shape-3D} shows three-dimensional views of two
solutions along the sixth branch. They have three full turns and nicely
illustrate the exact helical shape.

\begin{figure}
\begin{center}
{\bf~~~~~~~~(a)~~~~~~~~~~~~~~~~~~~~~~~~~~~~~~~~~~~~~~~~~~~~~~~~~~~(b)~~~}\\
\includegraphics[width=0.45\linewidth]{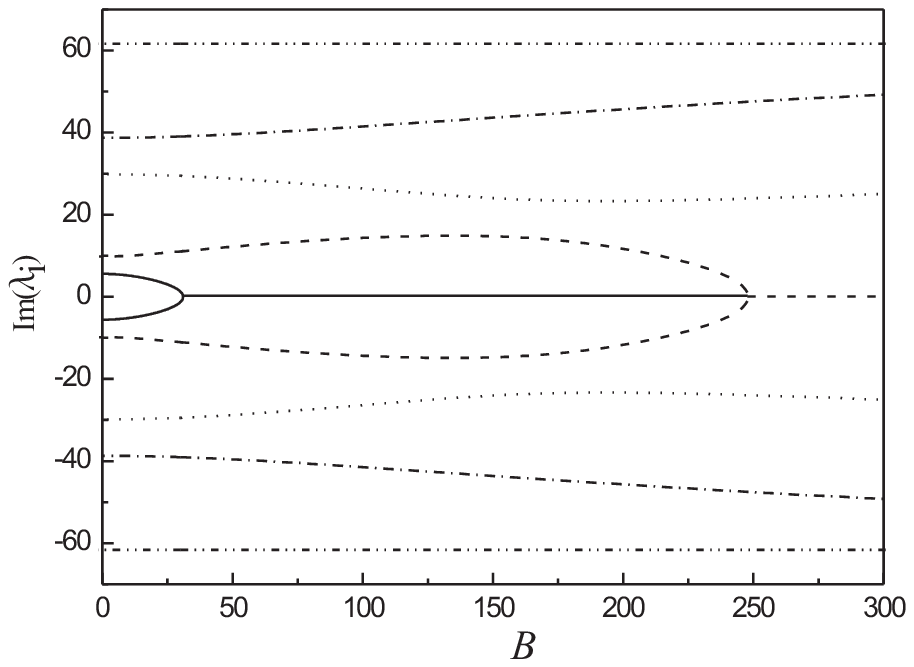}
~~~~~\includegraphics[width=0.45\linewidth]{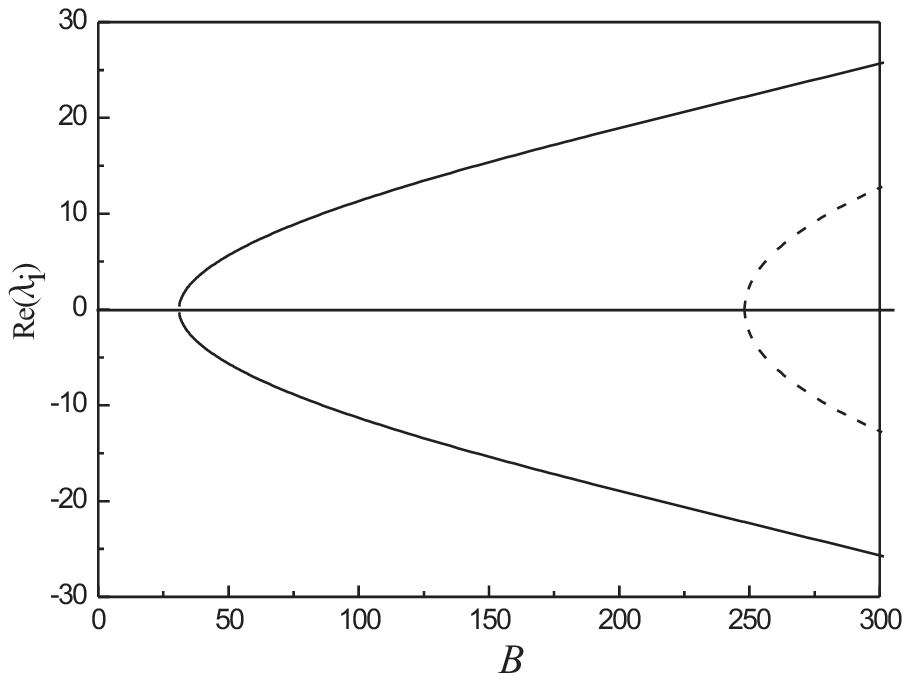}
\end{center}
\caption{Evolution of imaginary (a) and real (b) parts of the first five
pairs of eigenvalues along the trivial solution, from $B=0$ to $B>BP2$.}
\label{fig:hel:helix_w0-stab-0BP1BP2-imag}
\end{figure}

\begin{figure}
\begin{center}
\includegraphics[width=0.5\linewidth]{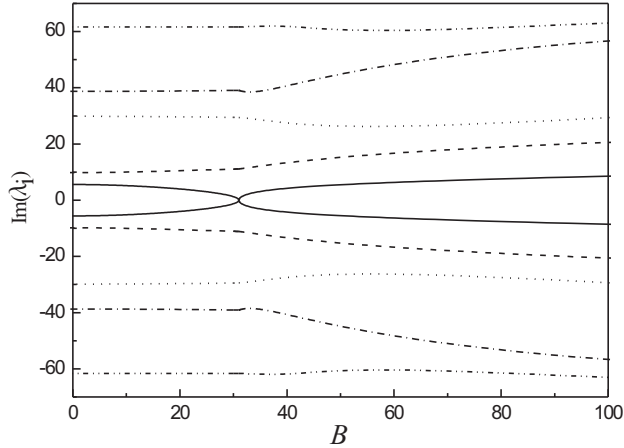}
\end{center}
\caption{Imaginary part evolution of the first five pairs of
eigenvalues along the trivial solution, from $B=0$ to $B=BP1$ and
then switching to branch b1. All real parts of the eigenvalues are zero.}
\label{fig:hel:helix_w0-stab-0BP1-b1-imag}
\end{figure}

\begin{figure}
\begin{center}
{\bf~~~~~~~(a)~~~~~~~~~~~~~~~~~~~~~~~~~~~~~~~~~~~~~~~~~~~~~~~~~(b)~~~}\\
\includegraphics[width=0.45\linewidth]{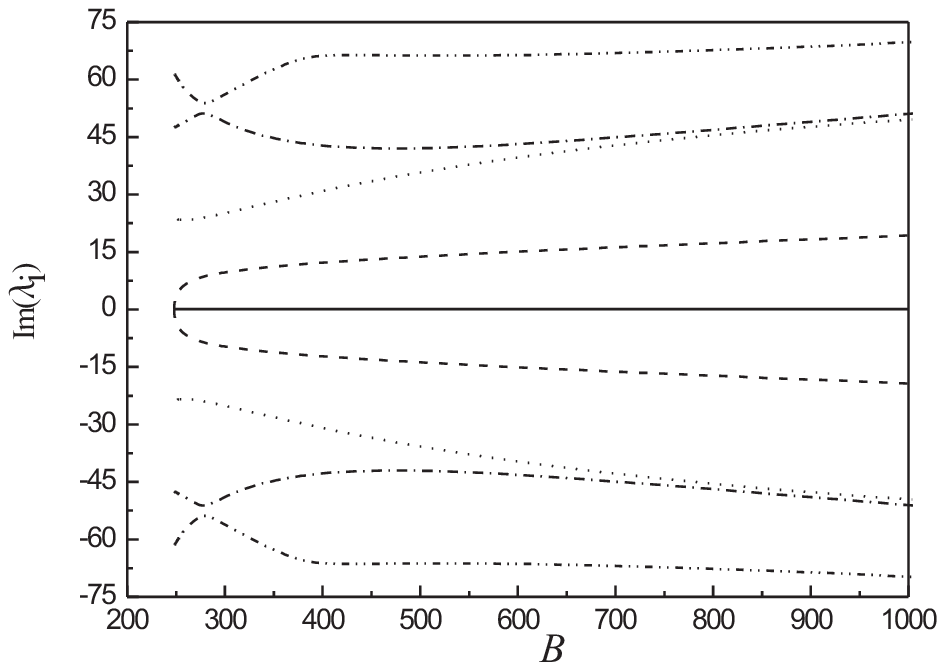}
~~~~~\includegraphics[width=0.45\linewidth]{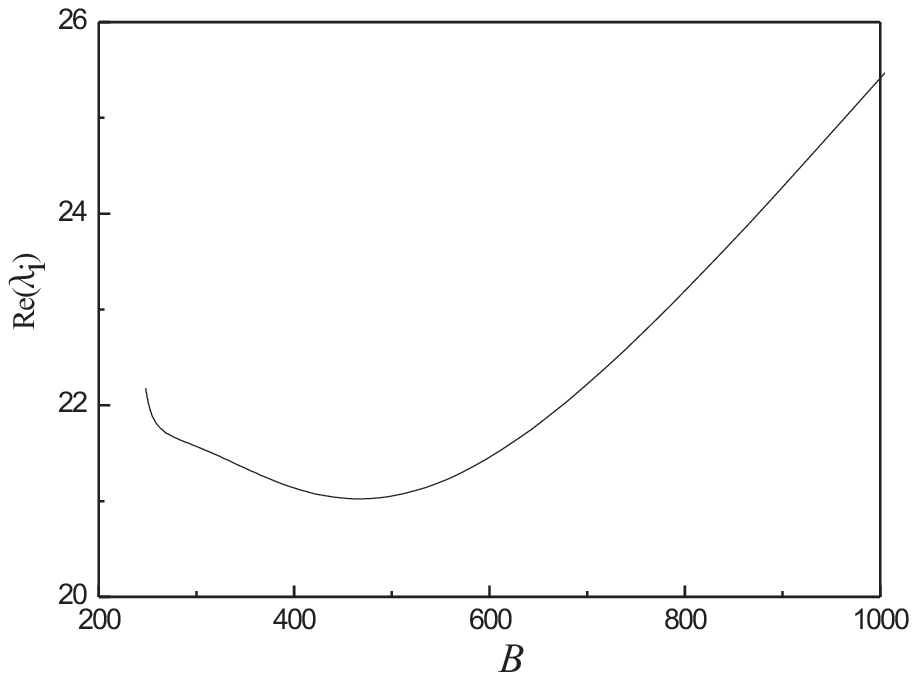}
\end{center}
\caption{Evolution of the imaginary (a) and (positive) real (b)
parts of the first five pairs of eigenvalues along branch b2.}
\label{fig:hel:helix_w0-stab-0BP2-b2-imag}
\end{figure}

Fig.~\ref{fig:hel:helix_w0-stab-0BP1BP2-imag} show the evolution of the
imaginary and real parts, respectively, of the first five pairs of
eigenvalues along the trivial solution, from $B=0$ to $B>BP2$. At successive
pitchfork bifurcations, pairs of imaginary eigenvalues collide at zero
and become real, one of the eigenvalues of the pair with positive real part,
signalling that the trivial solution becomes unstable at $B=BP1$, while a
further loss of stability occurs at $BP2$ where a second pair of eigenvalues
becomes real.

Fig.~\ref{fig:hel:helix_w0-stab-0BP1-b1-imag} shows the evolution
of the first five pairs of imaginary eigenvalues when switching at
BP1 from the trivial branch to b1, which is found to be stable. As
expected, solutions along b2 are found to be unstable (see
Figs~\ref{fig:hel:helix_w0-stab-0BP2-b2-imag}a,b).

\subsection{Whirling solutions ($\omega\neq0$) -- Hopf bifurcations}

\begin{figure}[h]
\begin{center}
\includegraphics[width=0.45\linewidth]{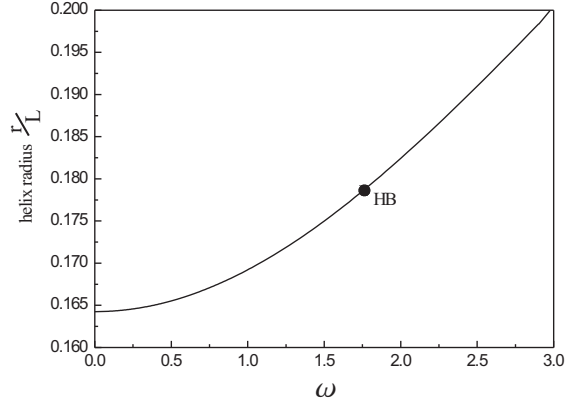}
\end{center}
\caption{Effect of $\omega$ on the helical radius of the solution
along b1 in Fig.~\ref{fig:hel:bif-helix} at $B=35$. The dot
indicates a Hopf bifurcation where the real parts of two complex
conjugate eigenvalues become positive (see
Fig.~\ref{fig:hel:helix_stab-b1_B35_gamma05_w-real}(a)) and the
solution loses stability.}
\label{fig:hel:bif_w_b35}
\end{figure}

\begin{figure}[b]
\begin{center}
\includegraphics[width=0.6\linewidth]{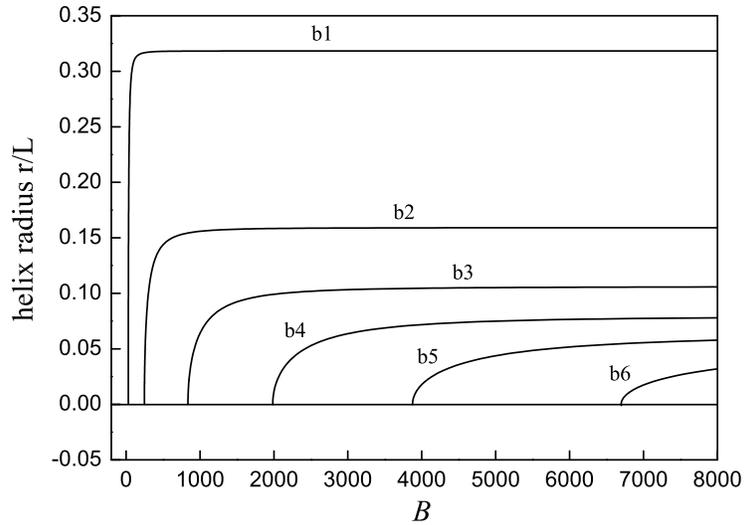}
\end{center}
\caption{Bifurcation diagram for a whirling isotropic rod subject to coat
hanger boundary conditions ($\omega=2$).}
\label{fig:hel:bif-diag-w2}
\end{figure}

\begin{figure}
\begin{center}
{\bf~~~~~~(a)~~~~~~~~~~~~~~~~~~~~~~~~~~~~~~~~~~~~~~~~~~~~~~~~~~(b)~~~}\\
\includegraphics[width=0.45\linewidth]{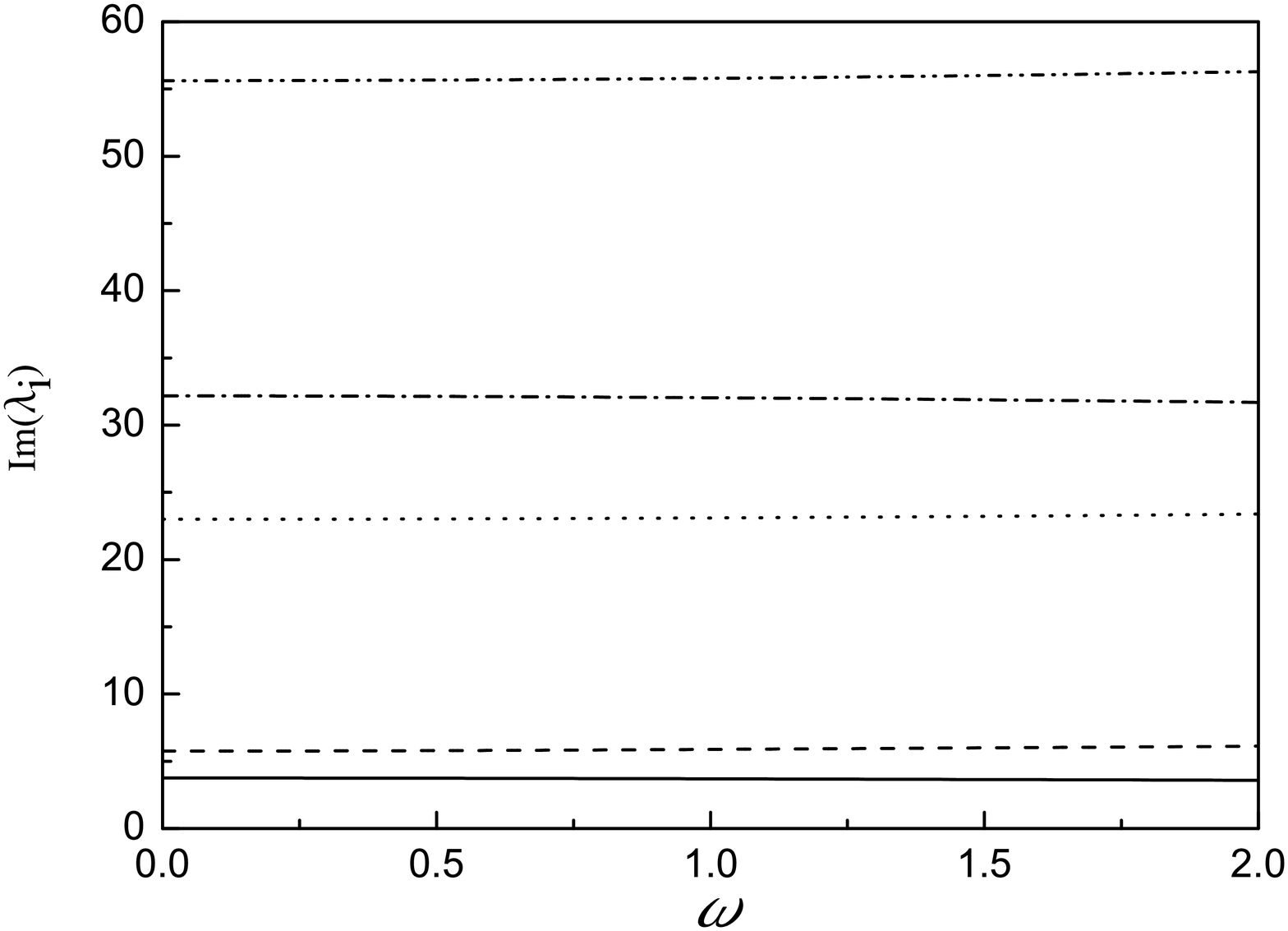}
~~~~~\includegraphics[width=0.45\linewidth]{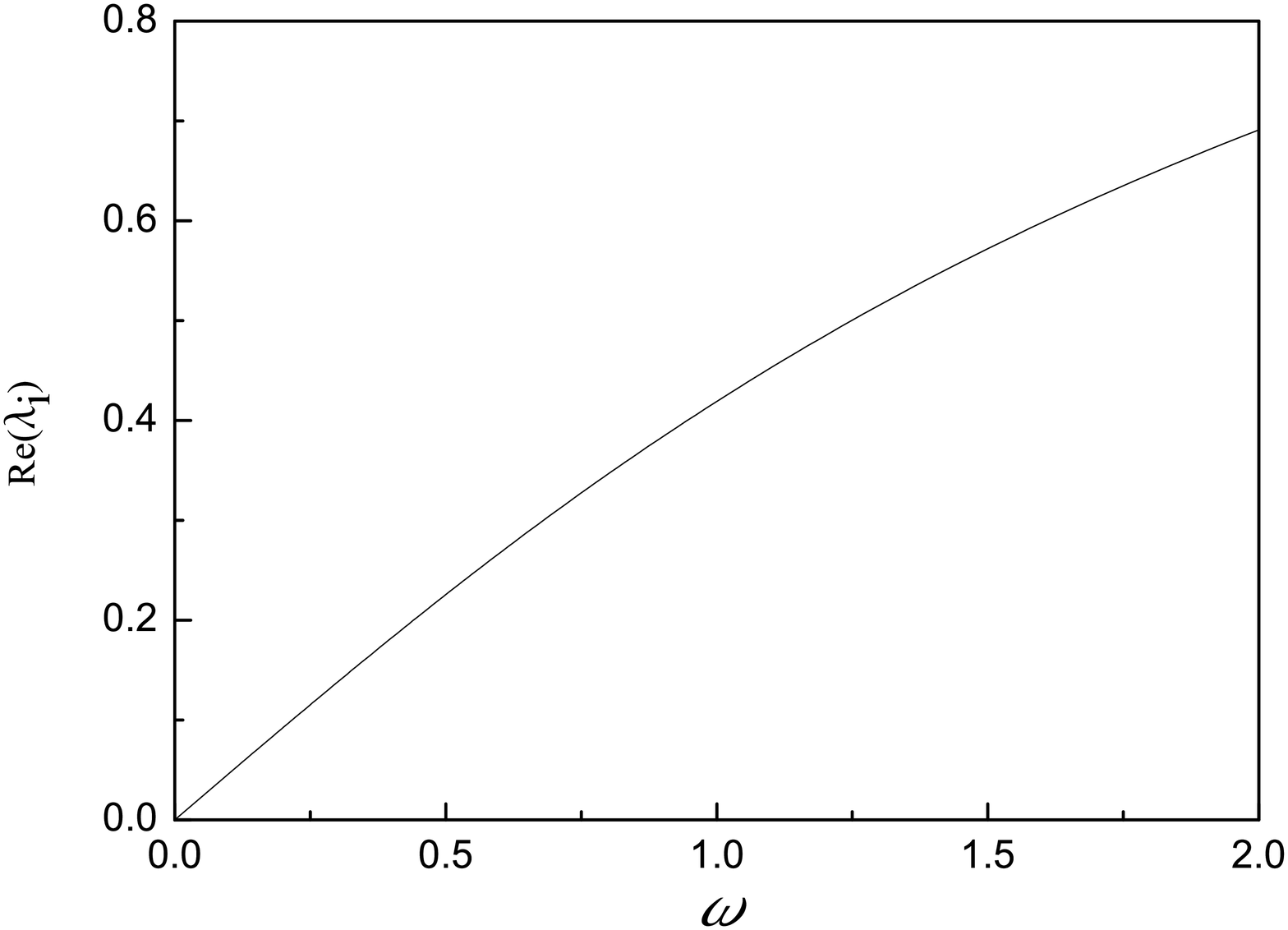} \vspace{0.2cm}\\
{\bf~~~~~~(c)~~~~~~~~~~~~~~~~~~~~~~~~~~~~~~~~~~~~~~~~~~~~~~~~~~(d)~~~} \vspace{0.1cm}\\
\includegraphics[width=0.45\linewidth]{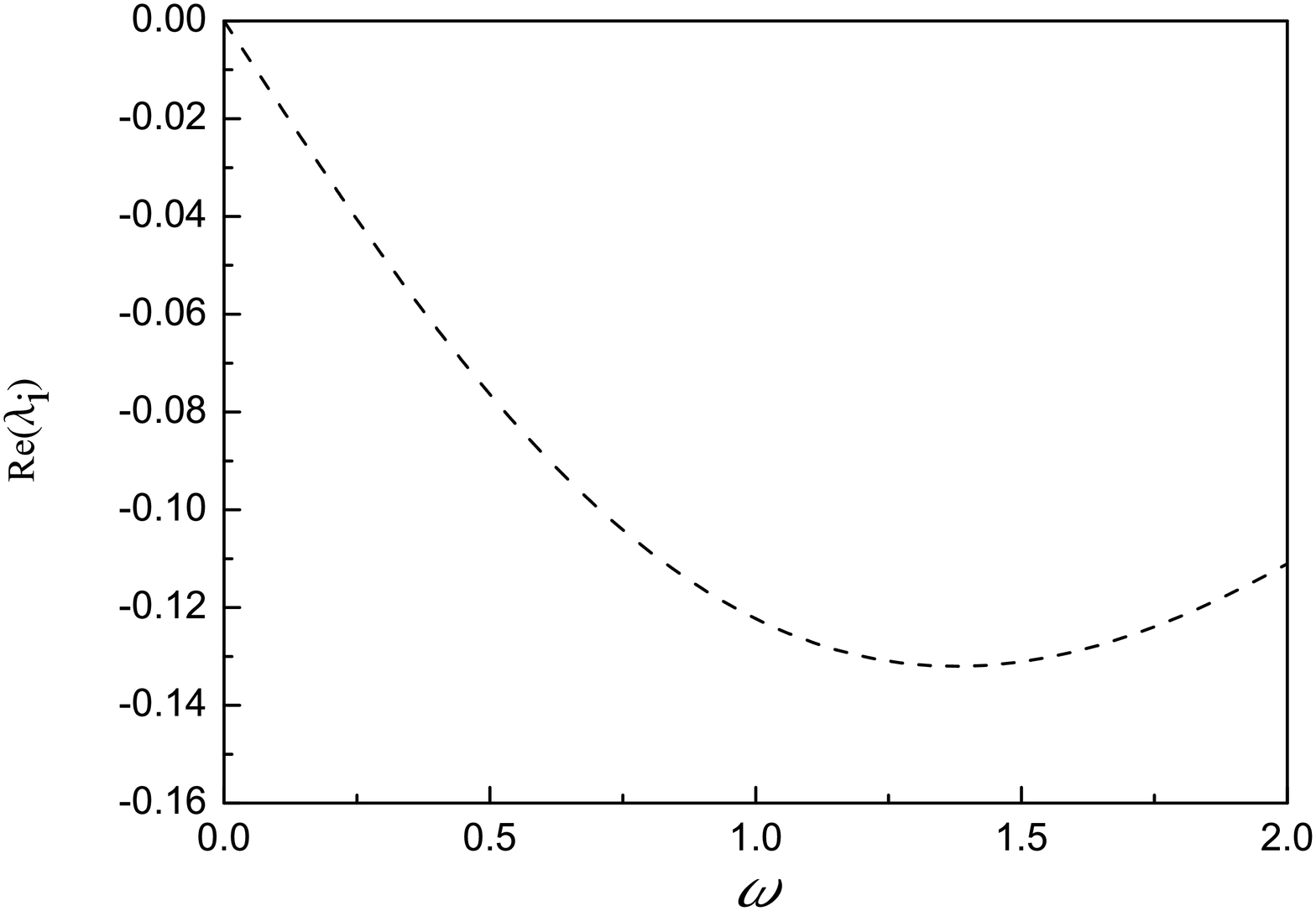}
~~~~~\includegraphics[width=0.45\linewidth]{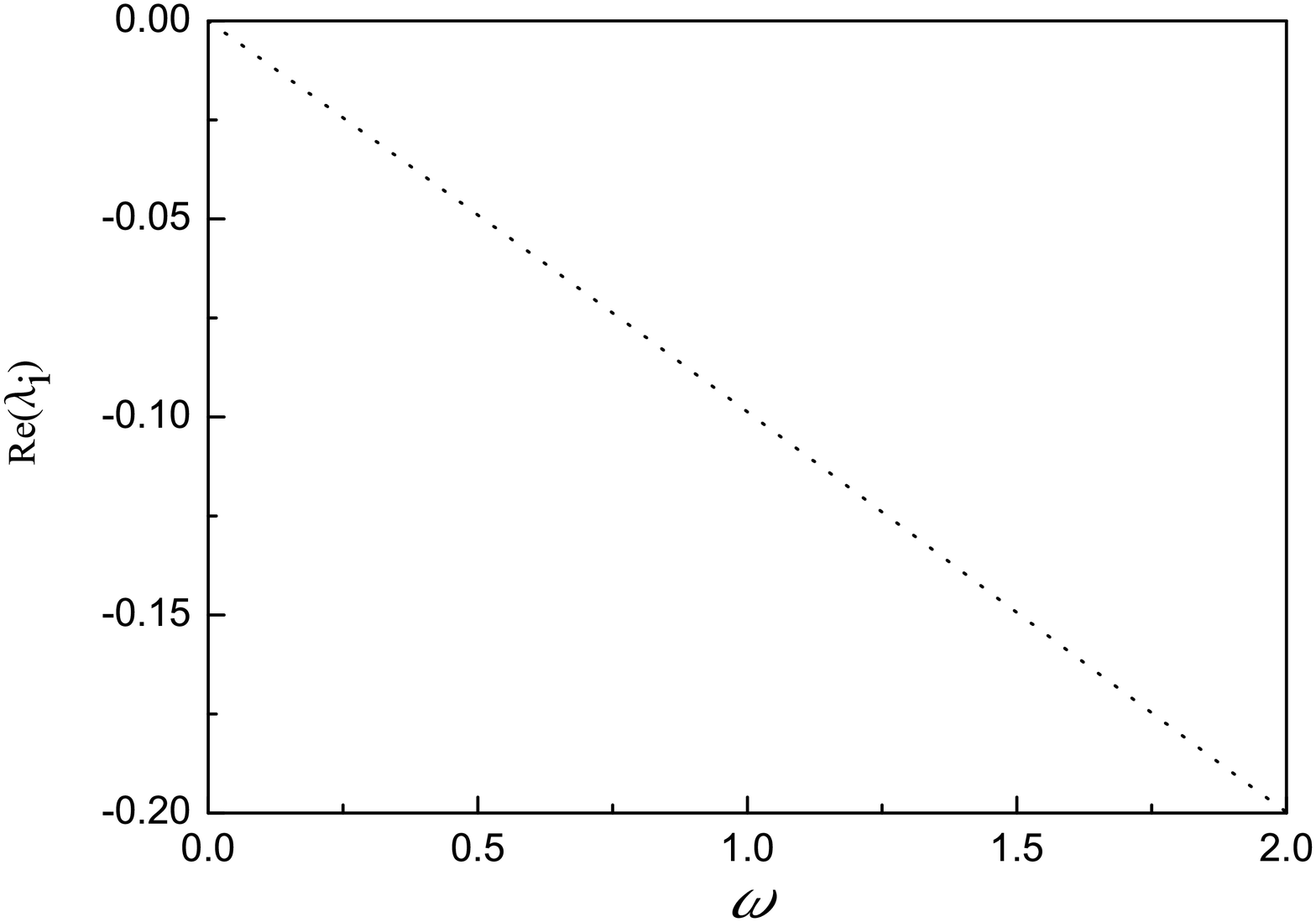} \vspace{0.2cm}\\
{\bf~~~~~~(e)~~~~~~~~~~~~~~~~~~~~~~~~~~~~~~~~~~~~~~~~~~~~~~~~~~(f)~~~} \vspace{0.1cm}\\
\includegraphics[width=0.45\linewidth]{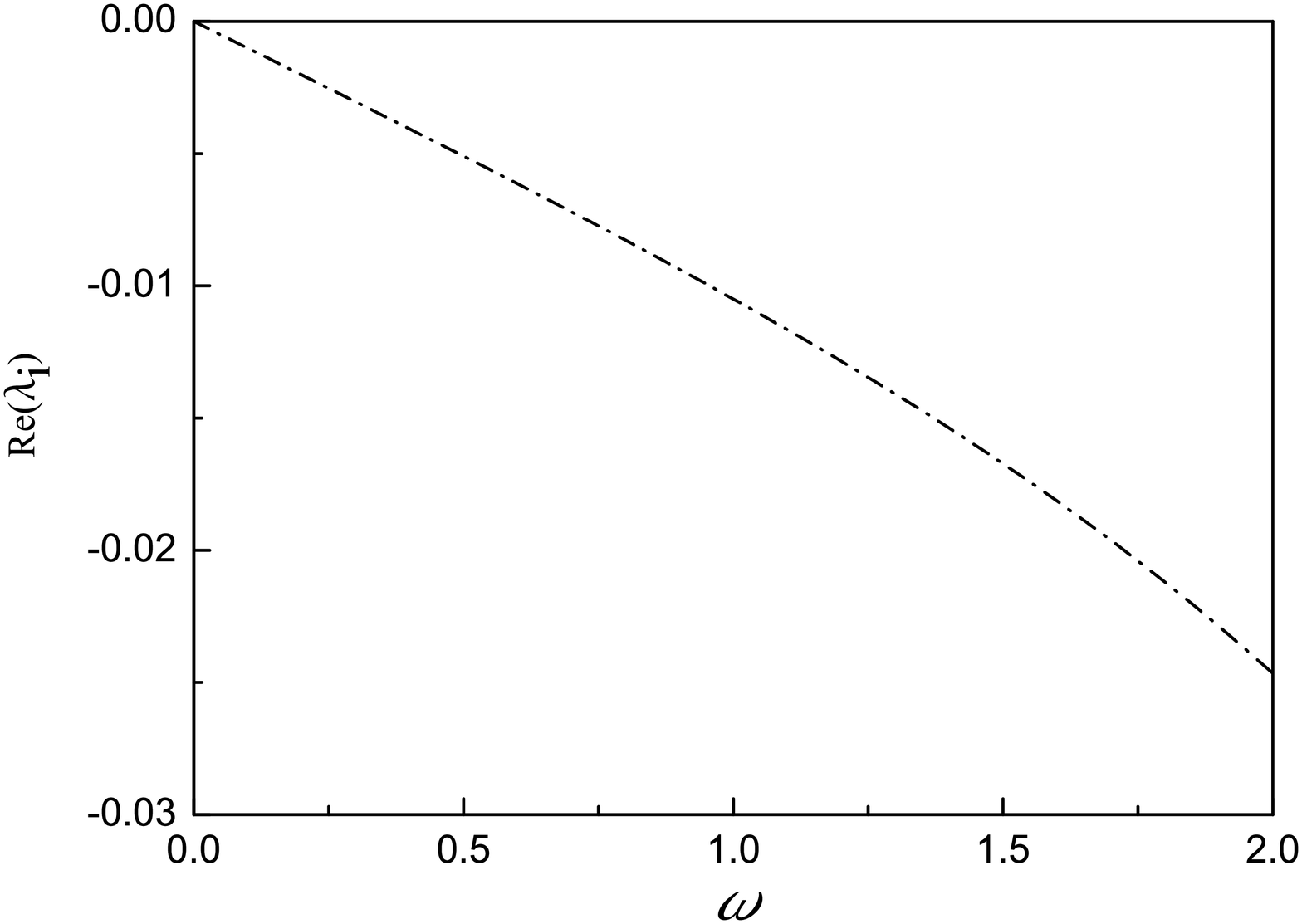}
~~~~~\includegraphics[width=0.45\linewidth]{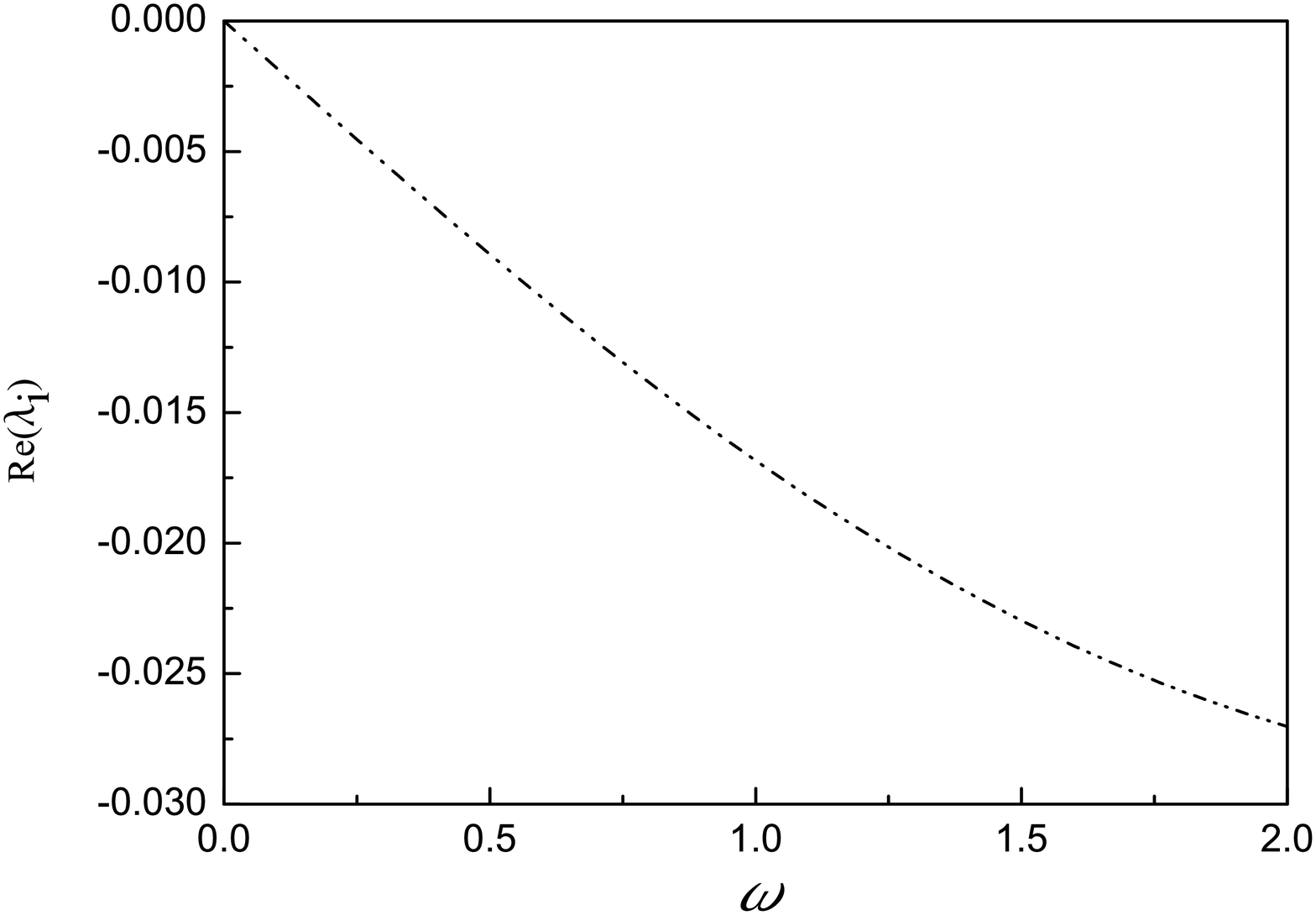}
\end{center}
\caption{Evolution of imaginary (a) and real part (b-f) of the first five
pairs of eigenvalues with respect to $\omega$ for a solution of the first
branch b1 of Fig.~\ref{fig:hel:bif-helix} for $B=35$, bifurcation diagram
in Fig.~\ref{fig:hel:bif_w_b35}. Note that first panel only
represents the positive imaginary eigenvalue, as the negative is
symmetric with respect to the real axis.}
\label{fig:hel:helix_w-stab-b1_real_imag}
\end{figure}

\begin{figure}
\begin{center}
{\bf~~~~~~(a)~~~~~~~~~~~~~~~~~~~~~~~~~~~~~~~~~~~~~~~~~~~~~~~~~~(b)~~~}\\
\includegraphics[width=0.45\linewidth]{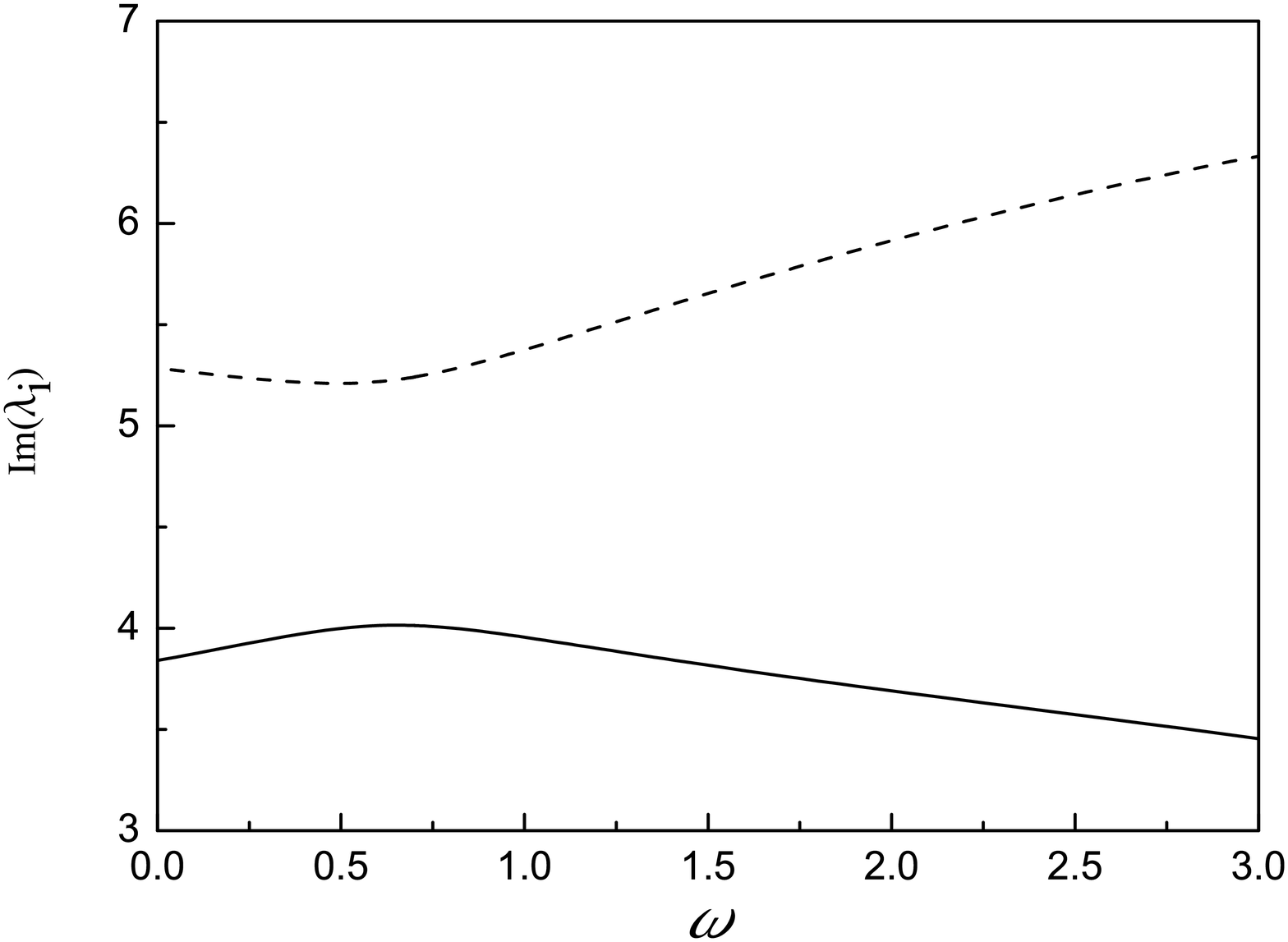}
~~~~~\includegraphics[width=0.45\linewidth]{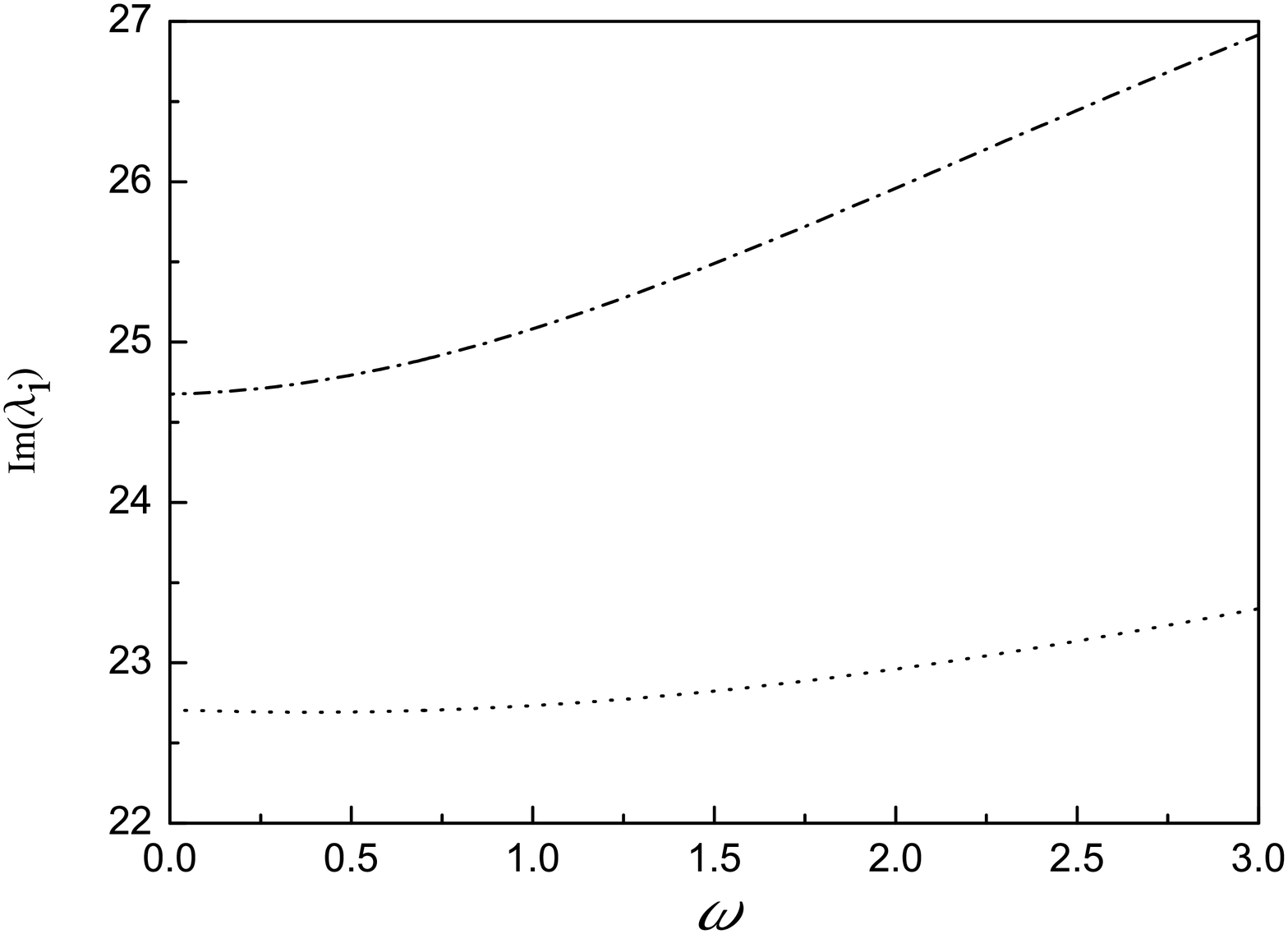} \vspace{0.2cm}\\
{\bf ~~~(c)} \vspace{0.1cm}\\
\includegraphics[width=0.45\linewidth]{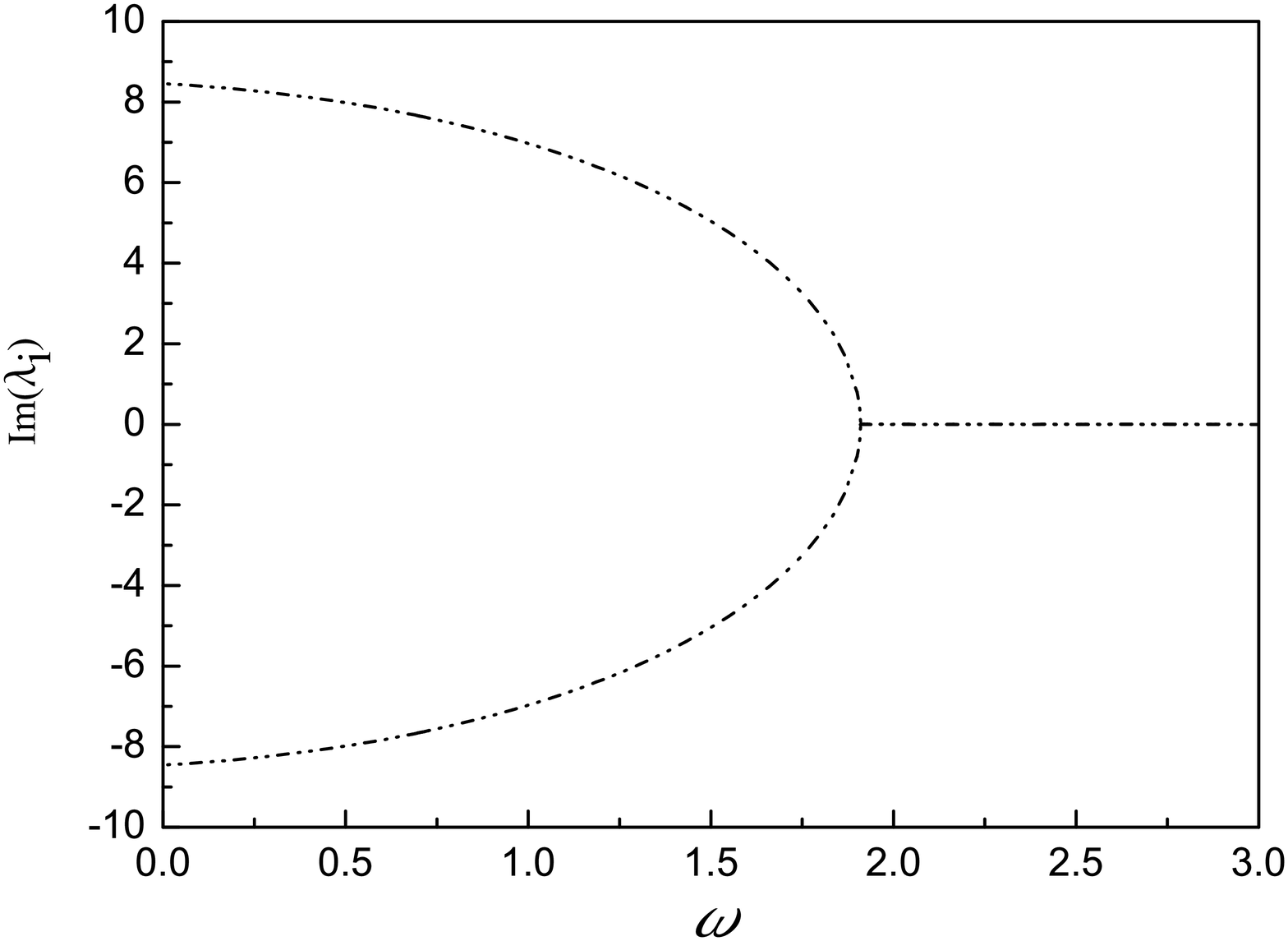}
\end{center}
\caption{Evolution with respect to $\omega$ of the imaginary part of
the first five pairs of eigenvalues for a solution along b1 at $B=35$ and
$\gamma=0.05$. The first two panels only show the positive imaginary part.}
\label{fig:hel:helix_stab-b1_B35_gamma05_w-imag}
\end{figure}

\begin{figure}
\begin{center}
{\bf~~~~~~(a)~~~~~~~~~~~~~~~~~~~~~~~~~~~~~~~~~~~~~~~~~~~~~~~~~~(b)~~~}\\
\includegraphics[width=0.45\linewidth]{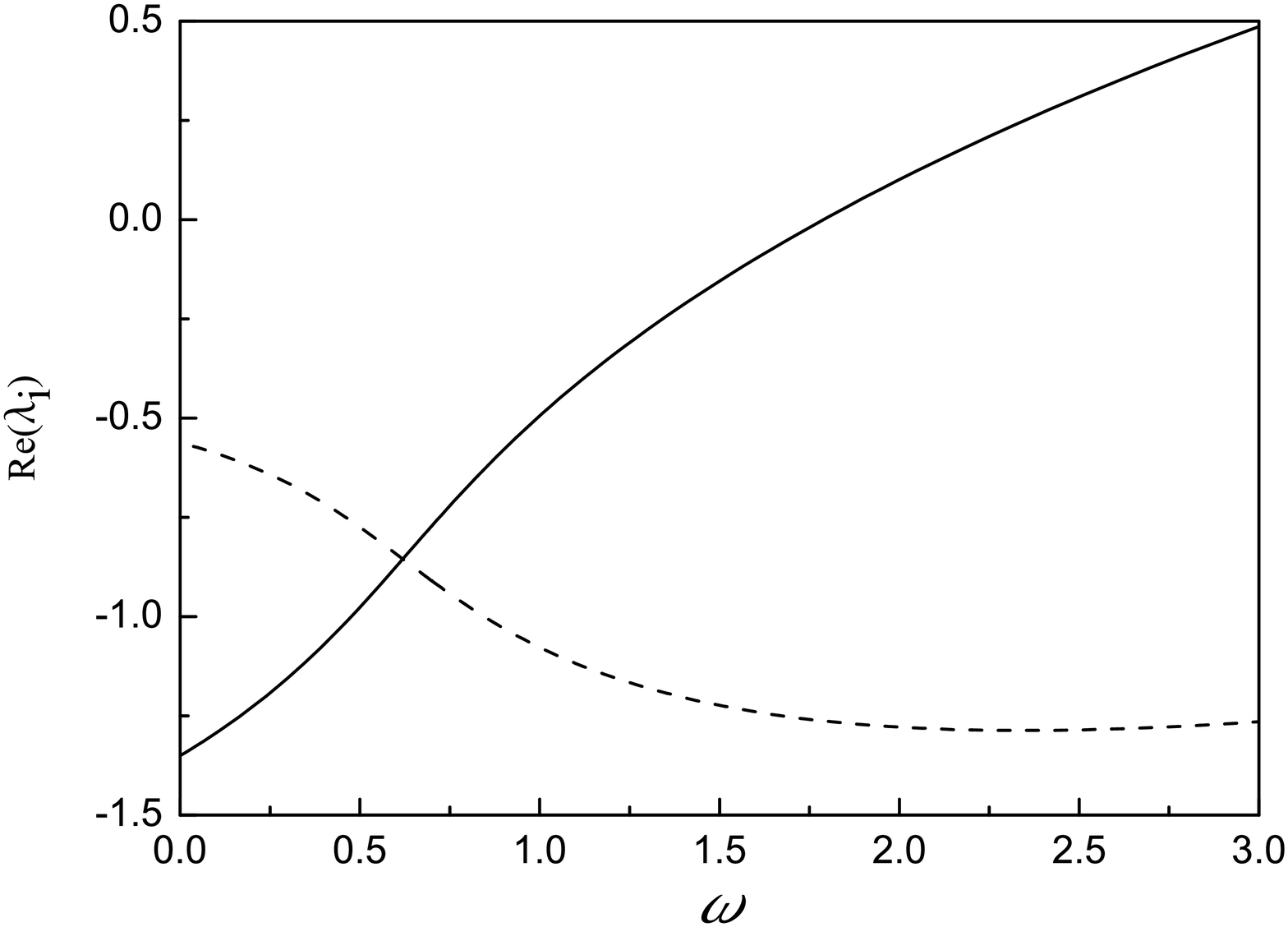}
~~~~~\includegraphics[width=0.45\linewidth]{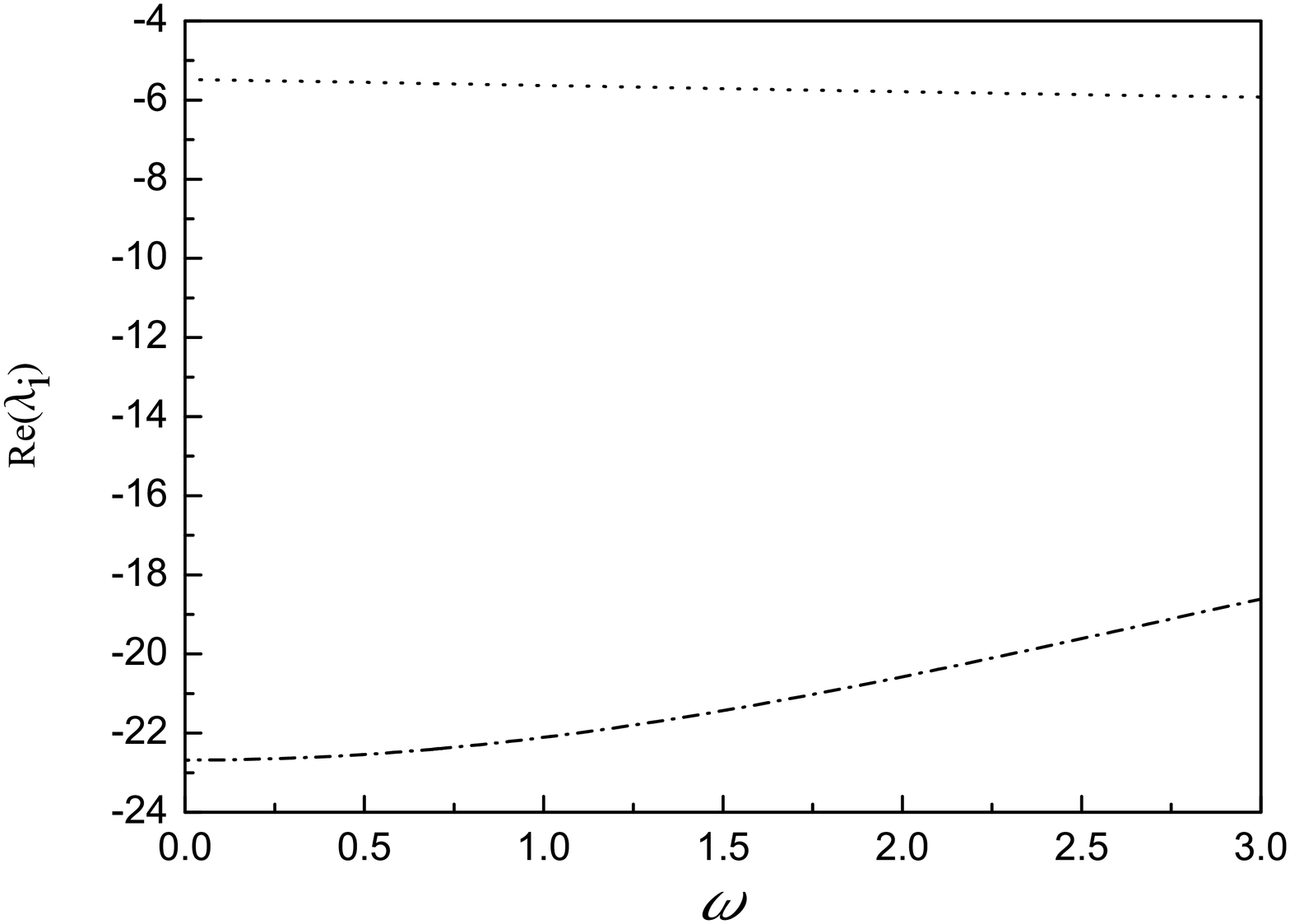} \vspace{0.2cm}\\
{\bf ~~~(c)} \vspace{0.1cm}\\
\includegraphics[width=0.45\linewidth]{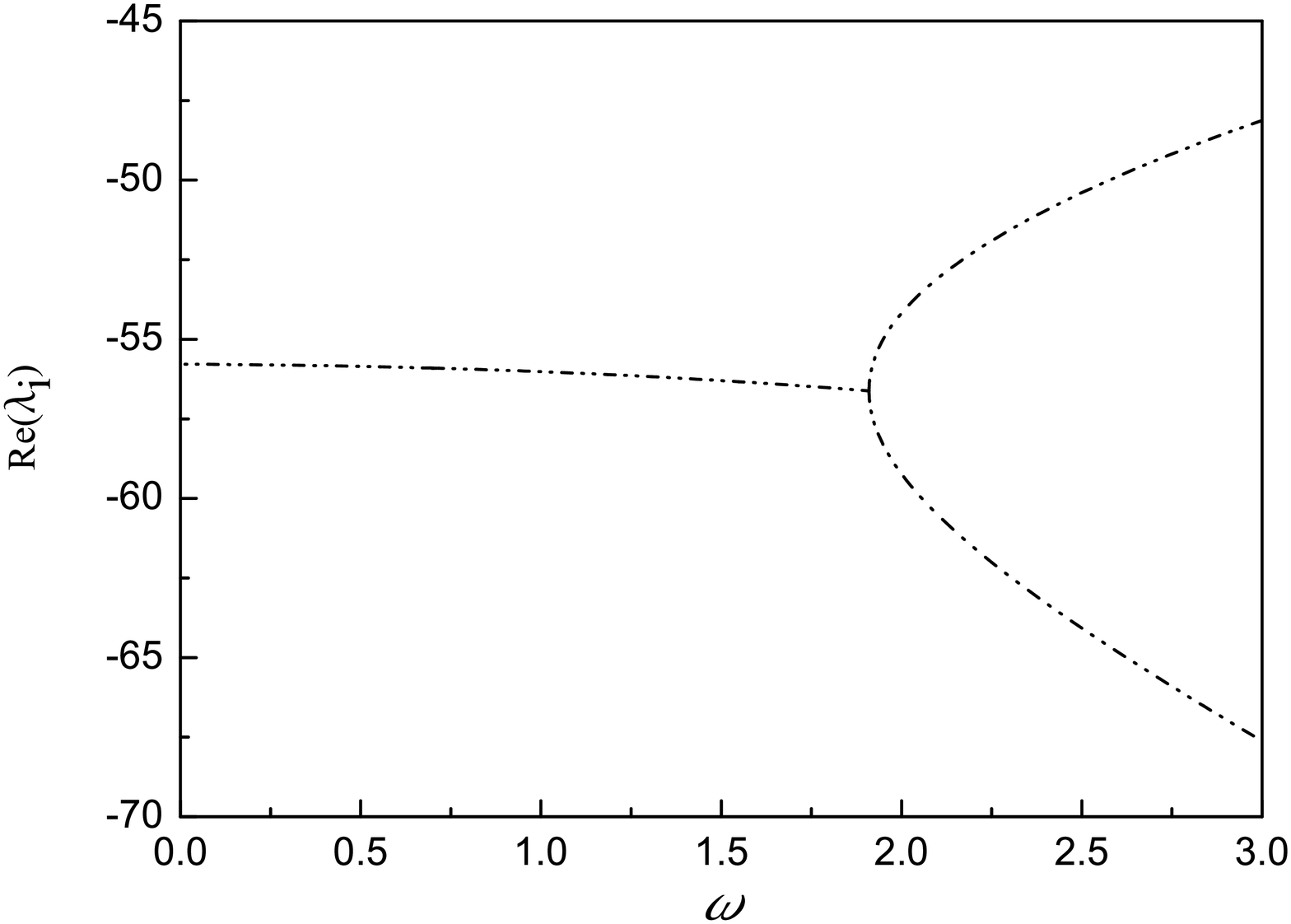}
\end{center}
\caption{Evolution with respect to $\omega$ of the real part of the
first five pairs of eigenvalues for a solution along b1 at $B=35$ and
$\gamma=0.05$.}
\label{fig:hel:helix_stab-b1_B35_gamma05_w-real}
\end{figure}

Next we consider whirling motions in which the supports $\bo v_0$
and $\bo v_1$ are spun about $\bo k$ with constant angular velocity
$\omega$. Solutions remain helical, as expected.
Fig.~\ref{fig:hel:bif_w_b35} shows the effect of $\omega$ on a
solution taken on the first bifurcating branch (b1) of
Fig.~\ref{fig:hel:bif-helix}. The helical radius increases with
$\omega$. Fig.~\ref{fig:hel:bif-diag-w2} shows the bifurcation
diagram for a fixed value of $\omega = 2$. Pitchfork bifurcations
along the trivial branch occur at $B=29.75$, 247.44, 836.80, etc.

\begin{figure}
\begin{center}
\includegraphics[width=0.5\linewidth]{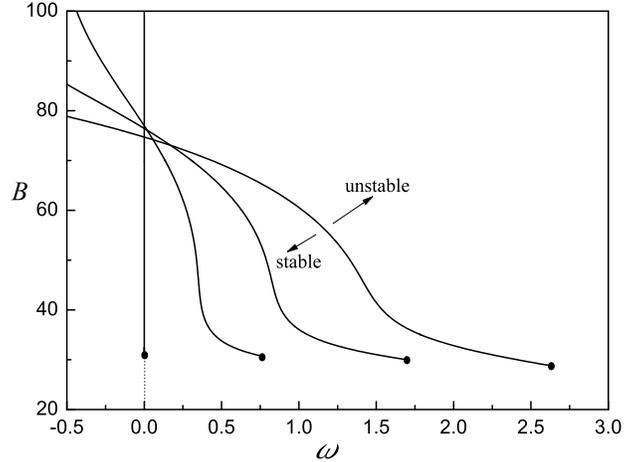}
\end{center}
\caption{Loci of Hopf bifurcations along b1 for different values of the
damping parameter $\gamma$: from left to right, $\gamma=0$, $\gamma=0.01$,
$\gamma=0.025$ and $\gamma=0.05$. Dots indicate points where the Hopf
bifurcation coalesces with the primary pitchfork bifurcation.}
\label{fig:hel:hopf-loci}
\end{figure}

The angular velocity tends to destabilise the helical solutions (see
Fig.~\ref{fig:hel:helix_w-stab-b1_real_imag}, where the eigenvalues are
shown as a function of $\omega$). We now introduce
damping by continuing in the parameter $\gamma$, and to further
investigate the stability of the first mode we fix
$\gamma=0.05$~\cite{valverde-cosserat} and perform continuation in
$\omega$. Figs~\ref{fig:hel:helix_stab-b1_B35_gamma05_w-imag} and
\ref{fig:hel:helix_stab-b1_B35_gamma05_w-real} show that at
$\omega=0$ the solution is stable and that around $\omega=1.75$ the
real part of the first eigenvalue becomes positive. Thus the system
loses stability in a Hopf bifurcation. The point where this occurs
is indicated in the bifurcation diagram in
Fig.~\ref{fig:hel:bif_w_b35}, which is still valid as $\gamma$ has
no effect on relative equilibria.
Fig.~\ref{fig:hel:hopf-loci} shows curves of Hopf bifurcations in the
$\omega$-$B$ parameter plane for various values of $\gamma$. At the
end points of these curves (indicated by dots) the Hopf bifurcation
coalesces with the pitchfork bifurcation in which the curve b1 is
created. Consequently, the whirling helical solution is stable below
these curves and unstable above. For small $\gamma$ the curve collapses
onto the line $\omega=0$. Note that for relatively large $\omega$
damping has a stabilising effect, but that for small positive and for
negative $\omega$ damping has a destabilising effect. The latter
behaviour is known from the classical linear stability theory of
gyroscopic systems \cite{ziegler}.

\subsection{The stationary anisotropic rod -- Hamiltonian-Hopf bifurcations}

\begin{figure}
\begin{center}
\includegraphics[width=0.6\linewidth]{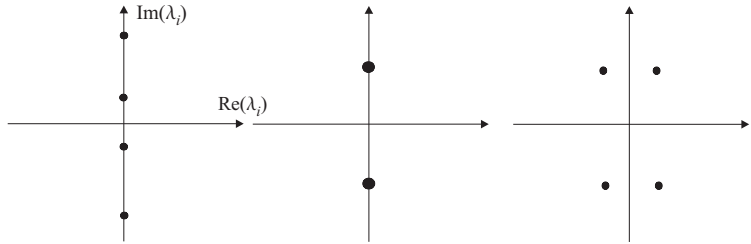}
\end{center}
\caption{Eigenvalue behaviour in a Hamiltonian-Hopf bifurcation.}
\label{fig:eigen-type}
\end{figure}

It is known that in the absence of damping and inertial effects
(i.e., $\gamma=0$, $\omega=0$) the equations for a rod in a magnetic
field have a Hamiltonian structure \cite{SH}. In Hamiltonian systems
a common mechanism for loss of stability is through a so-called
Hamiltonian-Hopf bifurcation \cite{van_der_meer}. In this
bifurcation two imaginary eigenvalues move on the imaginary axis,
meet at some non-zero value and then leave the axis to become a real
pair, as illustrated in Fig.~\ref{fig:eigen-type}. Since
eigenvalues come as conjugate pairs this event involves four
eigenvalues. In structural problems Hamiltonian-Hopf bifurcations
usually mark oscillatory instabilities such as flutter.

\begin{figure}
\begin{center}
\includegraphics[width=0.6\linewidth]{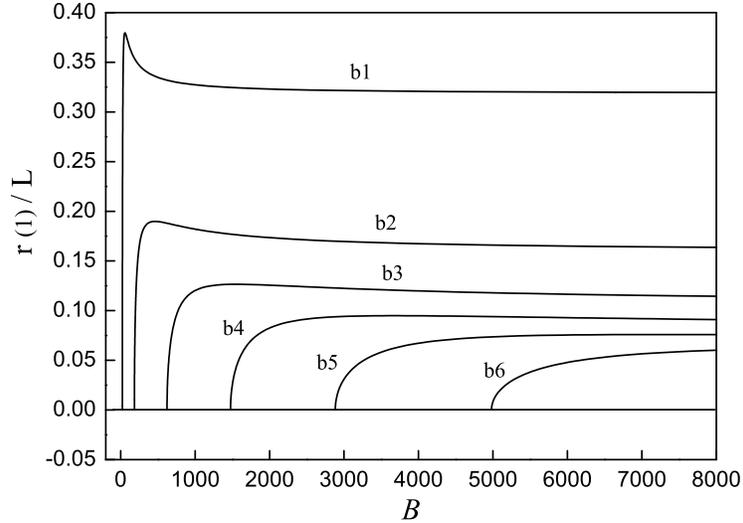}
\end{center}
\caption{Bifurcation diagram for a stationary anisotropic rod subject to
coat hanger boundary conditions ($\omega=0$, $R=0.5512$).
Bifurcating solutions are not exact helices.}
\label{fig:hel:bif-diag-R055w0}
\end{figure}

\begin{figure}
\begin{center}
{\bf~~~~~~(a)~~~~~~~~~~~~~~~~~~~~~~~~~~~~~~~~~~~~~~~~~~~~~~~~~~(b)~~~}\\
\includegraphics[width=0.45\linewidth]{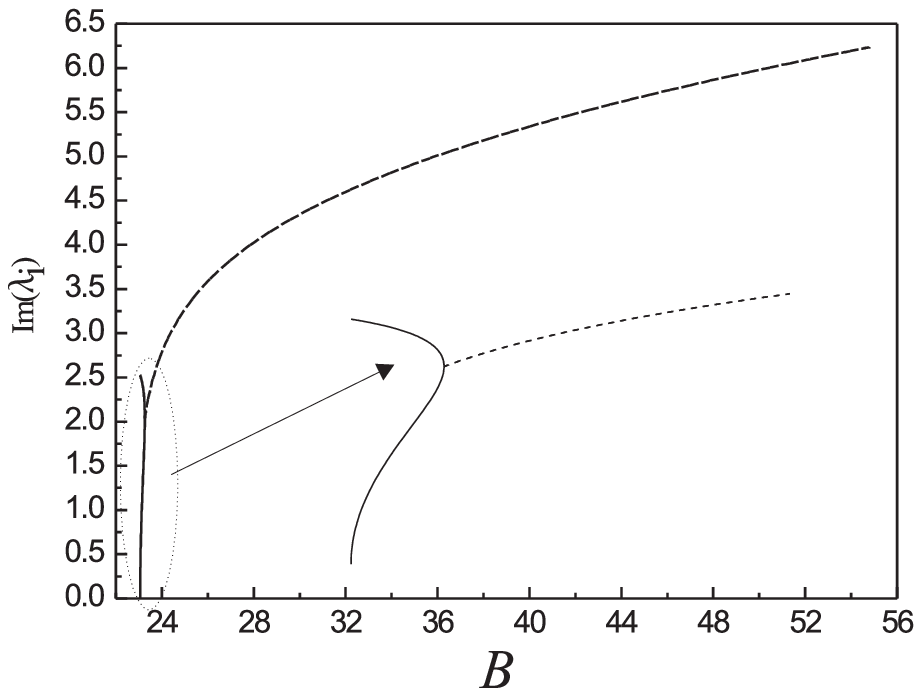}
~~~~~\includegraphics[width=0.45\linewidth]{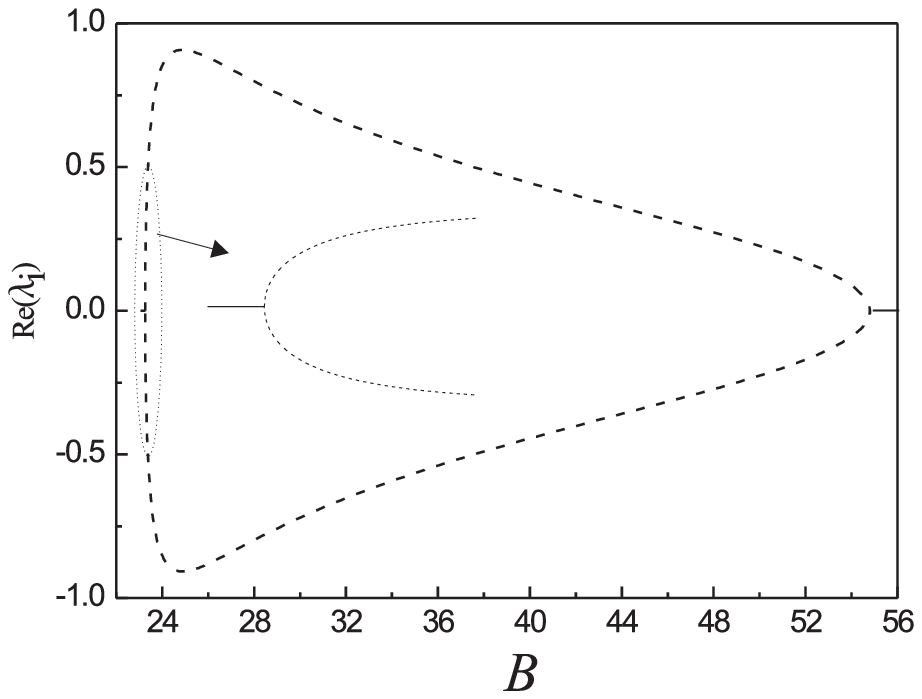}
\vspace{0.2cm}\\
{\bf~~~~~~(c)~~~~~~~~~~~~~~~~~~~~~~~~~~~~~~~~~~~~~~~~~~~~~~~~~~(d)~~~}\\
\includegraphics[width=0.45\linewidth]{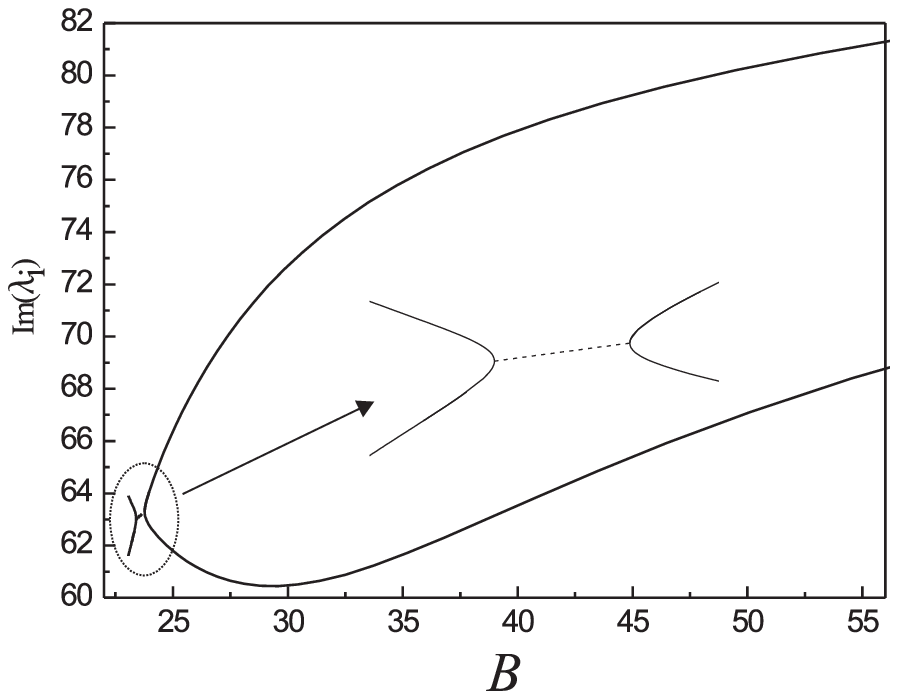}
~~~~~\includegraphics[width=0.45\linewidth]{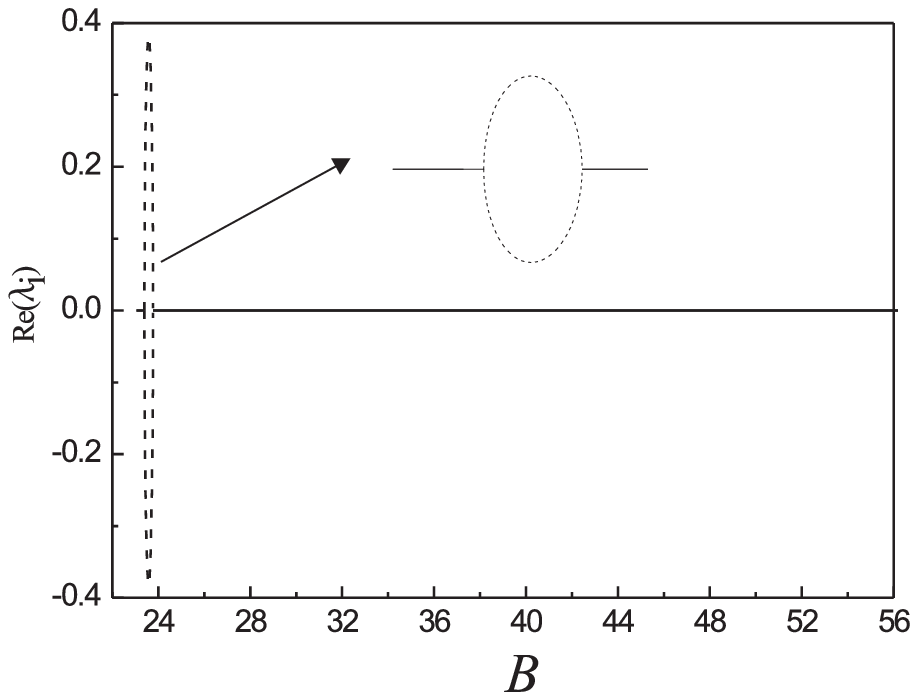}
\vspace{0.2cm}\\
\end{center}
\caption{Evolution with respect to $B$ of the first and second ((a)
and (b)) and fifth and sixth ((c) and (d)) eigenvalues for an
initially stable anisotropic-rod solution on the first bifurcating
branch (b1). Both pairs of eigenvalues undergo a Hamiltonian-Hopf
followed by a reverse Hamiltonian-Hopf bifurcation. ($\gamma=0$,
$\omega=0$, $R=0.55$.)}
\label{fig:hel:hopf1}
\end{figure}
\begin{figure}
\begin{center}
{\bf~~~~~~(a)~~~~~~~~~~~~~~~~~~~~~~~~~~~~~~~~~~~~~~~~~~~~~~~~~~(b)~~~}\\
\includegraphics[width=0.45\linewidth]{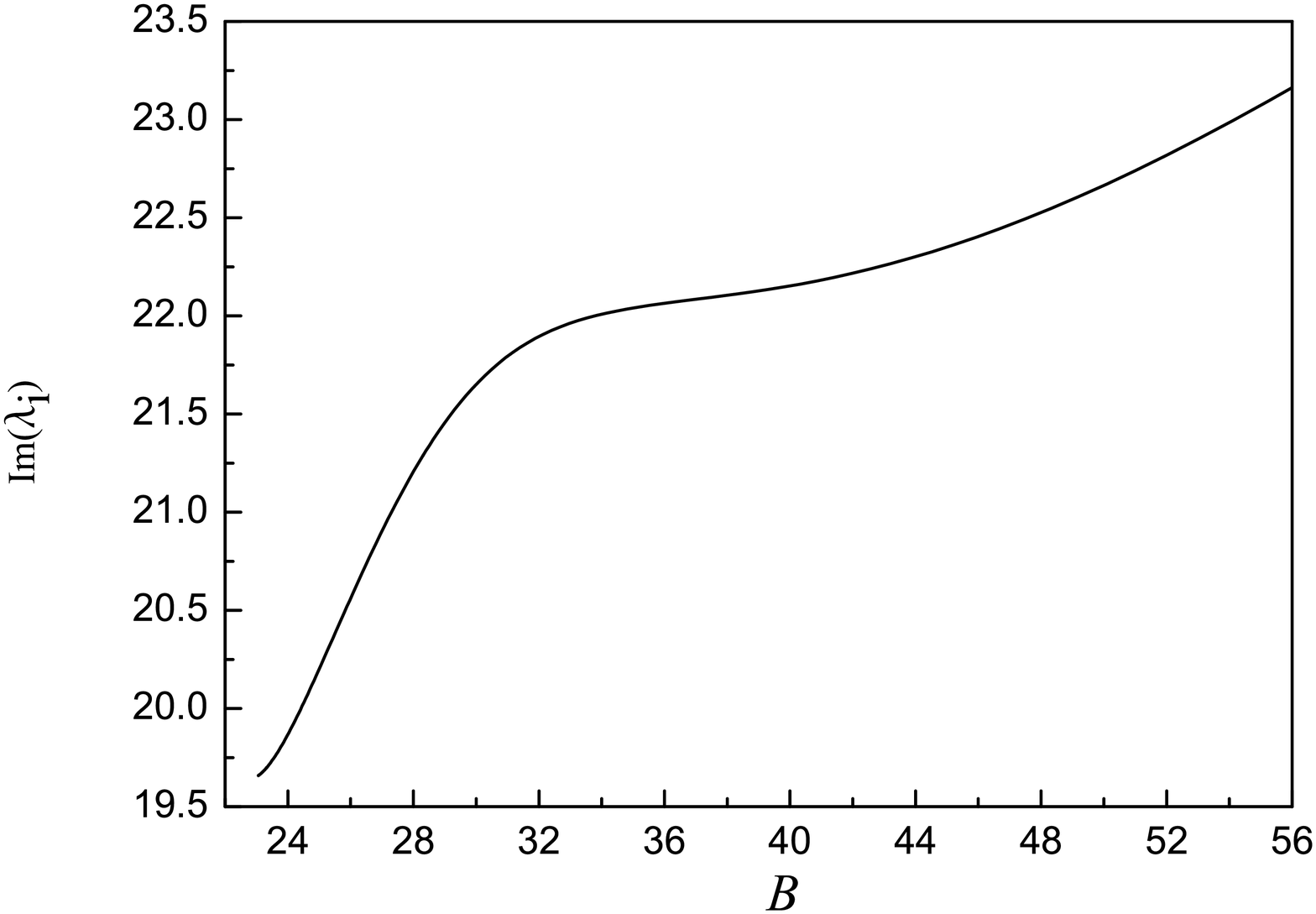}
~~~~~\includegraphics[width=0.45\linewidth]{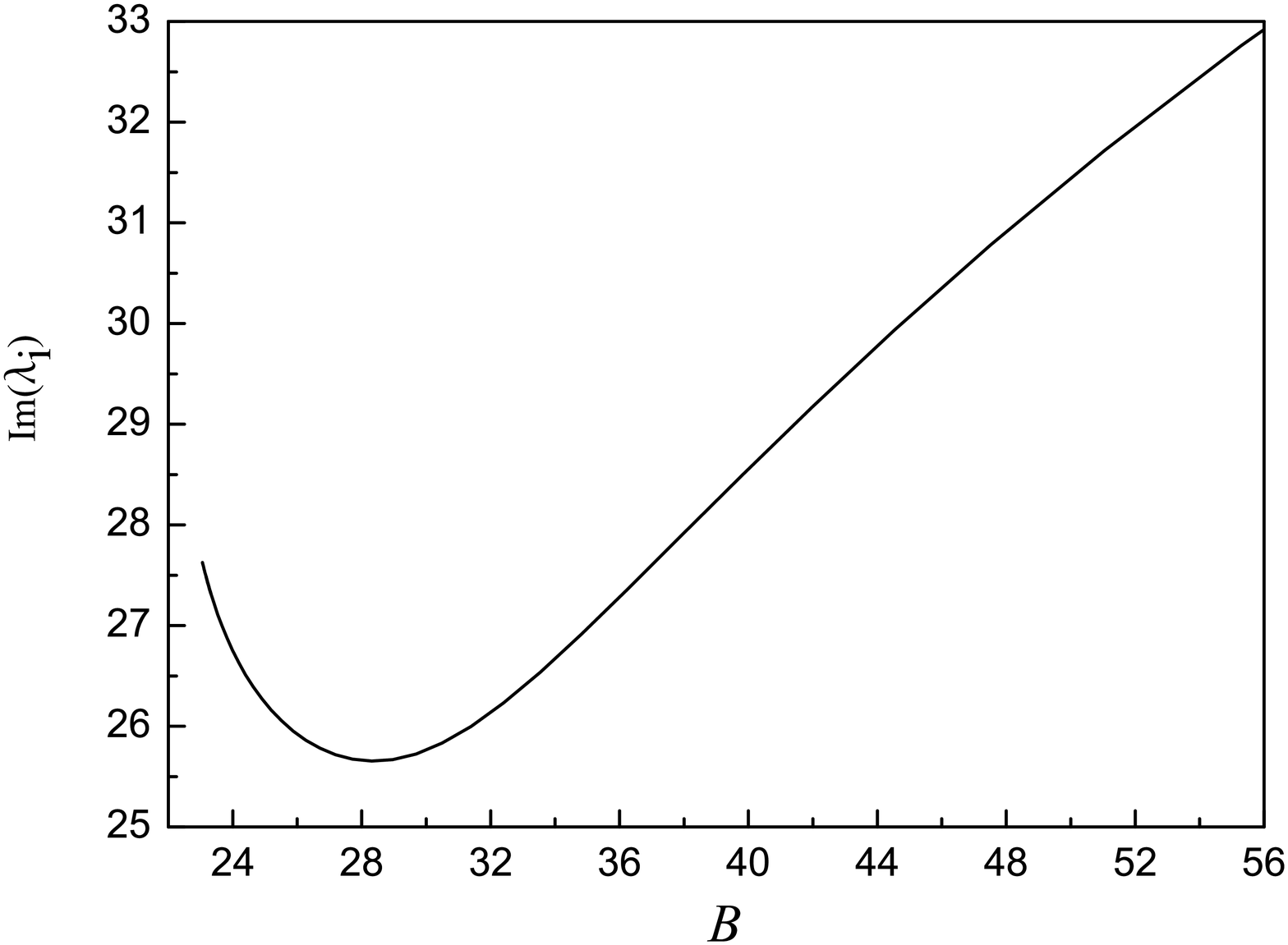}
\vspace{0.2cm}\\
{\bf ~~~(c)} \vspace{0.1cm}\\
\includegraphics[width=0.45\linewidth]{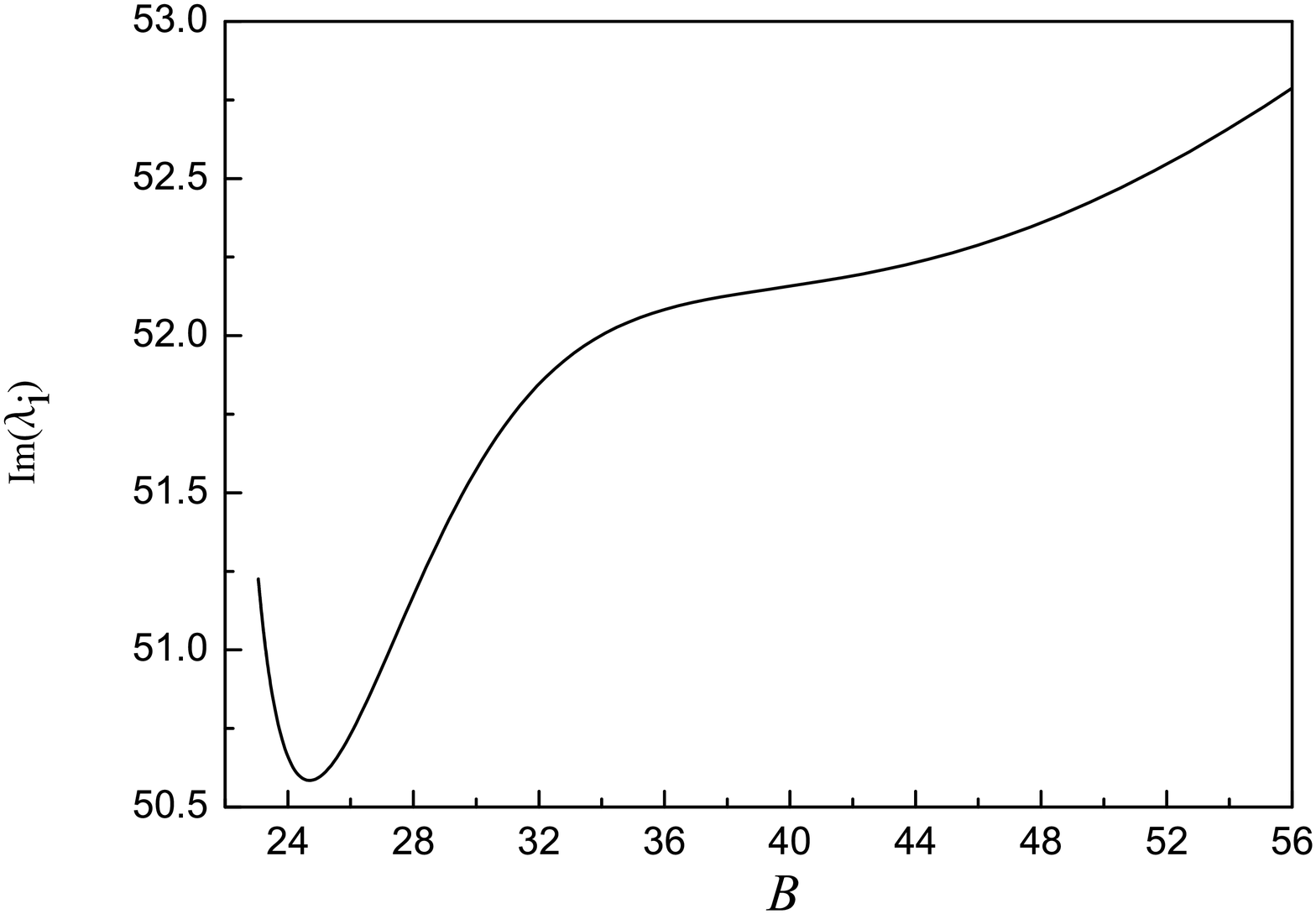}
\end{center}
\caption{Evolution with respect to $B$ ($\gamma=0$ and $\omega=0$)
of the imaginary part (real part equals zero) of the third, fourth
and seventh pairs of eigenvalues for non-circular
cross-section, $R=0.55$.}
\label{fig:hel:hopf2}
\end{figure}

\begin{figure}
\begin{center}
\includegraphics[width=0.45\linewidth]{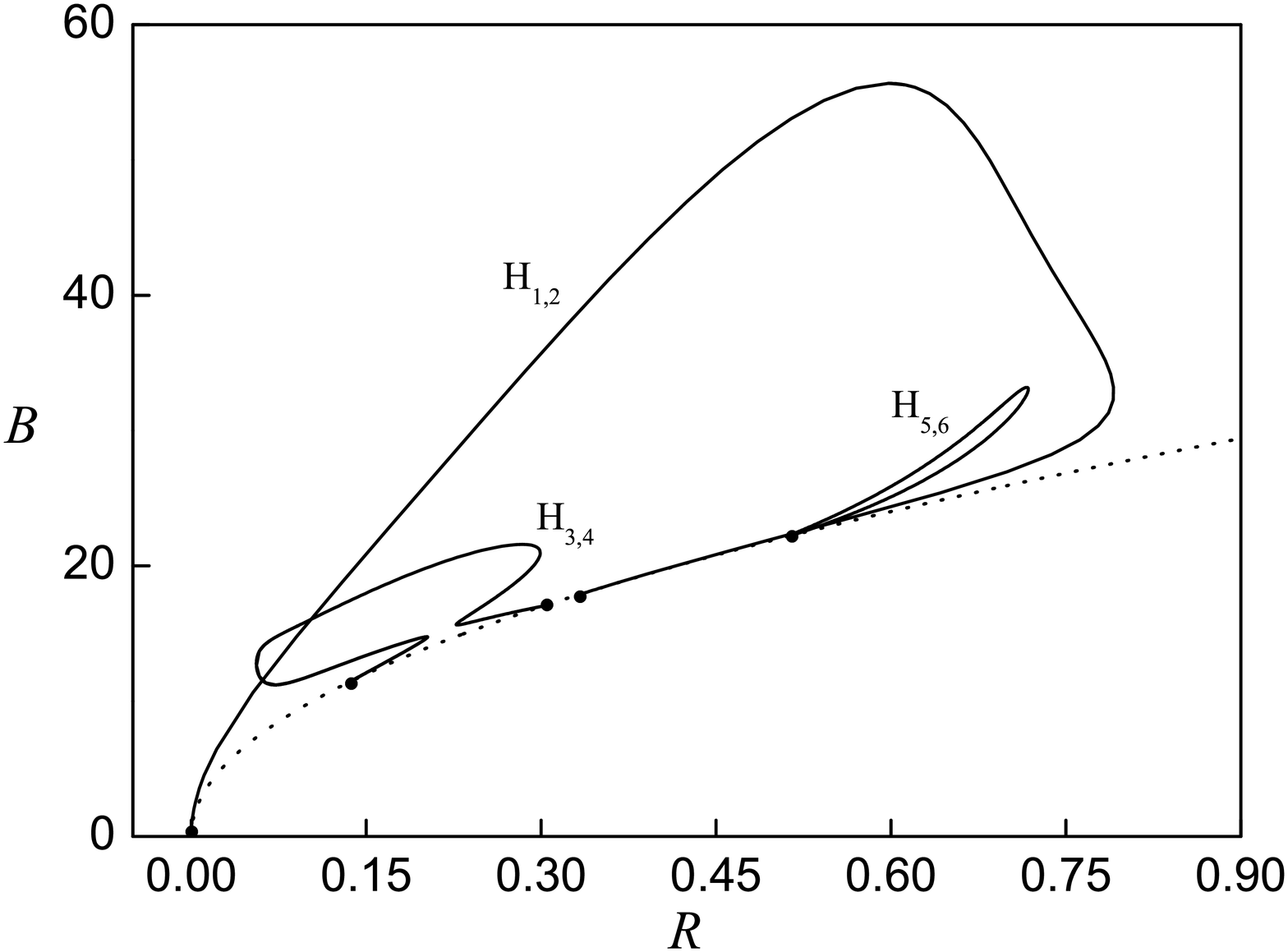}
\end{center}
\caption{Hamiltonian-Hopf loci for b1 solution at $\omega=0$, $R=0.55$. All
curves terminate on the dotted curve of pitchfork bifurcations given by
(\ref{char_coat}). We have instability of the solution inside the regions
bounded by solid and dotted lines. $H_{ij}$ labels the curve generated by
the $i$th and $j$th eigenvalue.}
\label{fig:hel:B_R-hopf}
\end{figure}

No Hamiltonian-Hopf bifurcation was found in our study of the isotropic
rod in Section \ref{coat_statics}. However, after breaking the
cross-sectional symmetry by taking $R<1$ Hamiltonian-Hopf
bifurcations are found to occur under $B$ continuation.
Fig.~\ref{fig:hel:bif-diag-R055w0} shows a bifurcation diagram for
$R=0.5512$. Pitchfork bifurcations along the trivial branch are found
at $B=23.058$, 184.56, 622.55, etc., in agreement with (\ref{char_coat}).
Bifurcating solutions are not pure helices.
Figs~\ref{fig:hel:hopf1} and \ref{fig:hel:hopf2} show
the eigenvalue behaviour for a solution on the first bifurcating
branch, b1, in Fig.~\ref{fig:hel:bif-diag-R055w0}. The first and
second eigenvalues collide at $B=23.26$, where the solution becomes
unstable. The fifth and sixth eigenvalues then collide at $B=23.39$
but shortly after become imaginary again in a reverse
Hamiltonian-Hopf bifurcation at $B=23.76$. The first and second
eigenvalues then also become imaginary again at $B=54.79$,
restabilising the solution (as far as we checked no other
eigenvalues cause instability, see Fig.~\ref{fig:hel:hopf2}, where
$\lambda_4$, $\lambda_5$ and $\lambda_7$ are shown).

Fig.~\ref{fig:hel:B_R-hopf} shows loci of Hamiltonian-Hopf bifurcations
in the $R$-$B$ parameter plane. At the end points of the curves the
Hamiltonian-Hopf bifurcation coalesces with the pitchfork bifurcation.
The curve of pitchfork bifurcations is given by (\ref{char_coat}) and
is included in dotted lines. We have instability inside the region
bounded by solid and dotted curves. Note that the stability of the b1
solution is almost entirely, but not completely, determined by the
first two eigenvalues.

\section{Conclusion}

We have shown that whirling current-carrying transversely isotropic
rods bifurcate under increasing magnetic field (or current) into
exact helical shapes provided one applies what we call coat hanger
boundary conditions. The first bifurcating branch, containing
solutions with half a helical turn, is stable while higher-order
branches are all unstable.

We stress that the stability analysis of Section 5.1 together with the
helical analysis in Appendix B gives a complete picture of
magnetically-induced helical buckling subject to these boundary
conditions. In Section 5.1 we show that the critical loads for a straight
rod are given by equation (\ref{char_coat}). These bifurcation points
are non-degenerate and valid for both isotropic and anisotropic rods.
The analysis says nothing, however, about the type of solutions that
bifurcate. They could be helices or not. In Appendix B we then compute
branches of helical solutions and show that, for isotropic rods (i.e.,
$R=1$), they intersect the trivial branch of straight rod solutions
precisely at the critical loads computed in Section 5.1. Taken together
these results prove that for an isotropic rod all solutions bifurcating
from the straight rod are helical. For an anisotropic rod, on the other
hand, the critical loads are still given by (\ref{char_coat}) (with
$R\neq 1$) but no branch of helical solutions in Appendix B is found to
intersect the trivial branch and therefore the bifurcating solutions are
non-helical, as confirmed numerically in Section 6.3.

We have also investigated the stability of post-buckling solutions and
found Hopf bifurcations where stable helical solutions lose stability.
In the case of a non-rotating anisotropic rod we found secondary
instabilities given by Hamiltonian-Hopf bifurcations. Unfortunately, our
method allows us only to study stationary or quasi-stationary (whirling)
solutions, so it is not clear what type of stable solutions occur after
these secondary bifurcations. To investigate this one would have to do
simulations based on direct discretisation of the PDEs
(\ref{eq:lmomentum4}) and (\ref{eq:amomentum4}).

For whirling isotropic rods subject to welded boundary conditions
magnetic buckling is described by a doubly-degenerate pitchfork
bifurcation (see \cite{clamped-JNLS}, where we also showed that the
same is true for non-rotating states of anisotropic rods). This is
because the equations are invariant under rotation about the axis
($\bo e_3$ axis) of the supports. This symmetry property complicated
Wolfe's analysis (for non-rotating states of isotropic rods), which
had to take account of the variational nature of the problem to prove
existence of non-trivial bifurcating states \cite{wolfe3}. By contrast,
the coat hanger boundary conditions here introduced break the $S^1$
symmetry down to $Z_2$ symmetry (reflection symmetry along $\bo v_0$
and $\bo v_1$), and no problems in the application of standard results
from bifurcation theory should arise. Indeed, we find the (isolated)
critical values of the magnetic field to be given by a remarkably
simple explicit expression.

Helical solutions are often used and studied in applications of elastic
rods or filaments. These solutions are incompatible with the common
(aligned) clamped, pinned, Cardan joint and other boundary conditions.
Boundary conditions are therefore often ignored in these studies, making
it impossible to do a stability analysis. Here we have introduced a set
of boundary conditions that does allow for helical solutions, not only
in magnetic buckling but also in traditional buckling due to compression,
twist or whirl. These boundary conditions allow one to study exact
helical solutions in finite-length rods. We have also proposed a
mechanical device that can be used for generating helical solutions in
the laboratory.

\section*{Acknowledgement}
J.V. wishes to thank the Andalusian Regional Ministry for financial
support through the Excellence Research Programme, under the FQM-4239
project.

%%%%%%%%%%%%%%%%%%%%%%%%%%%%%%%%%%%%%%%%%%%%%%%%%%%%%%%%%%%%%%%%%%%%%%
%%%%%%%%%%%%%%%%%%%%%%%%%%%%%%%%%%%%%%%%%%%%%%%%%%%%%%%%%%%%%%%%%%%%%%

%%%%%%%%%%%%%%%%%%%%%%%%%%%%%%%%%%%%%%%%%%%%%%%%%%%%%%%%%%%%%%%%%%%%%%
% Here's where you specify the bibliography style file.
% The full file name for the bibliography style file
% used for an ASME paper is asmems4.bst.
%\bibliographystyle{asmems4}

%%%%%%%%%%%%%%%%%%%%%%%%%%%%%%%%%%%%%%%%%%%%%%%%%%%%%%%%%%%%%%%%%%%%%%
% The bibliography is stored in an external database file
% in the BibTeX format (file_name.bib).  The bibliography is
% created by the following command and it will appear in this
% position in the document. You may, of course, create your
% own bibliography by using thebibliography environment as in
%

\newpage

\appendix
\begin{center}
\textbf{Appendix A: Matrices for the linearisation}
\end{center}

\noindent
The matrices $\bo B_i$ appearing in equation (\ref{per:lmomentum}) are given
by

\begin{equation}
\bo B_1 = \left(\begin{array}{ccc}
 0& -\kappa_3^0 & \kappa_2^0\\
 \kappa_3^0& 0 & -\kappa_1^0\\
 -\kappa_2^0& \kappa_1^0 &0
\end{array}\right),
\nonumber
\end{equation}

\begin{equation}
\bo B_2 = \omega^2 \left(\begin{array}{ccc}
 d_{11}^0& d_{12}^0 & 0\\
 d_{21}^0& d_{22}^0 & 0\\
 d_{31}^0& d_{32}^0 & 0
\end{array}\right),
\nonumber
\end{equation}

\begin{equation}
\bo B_3 =  \left(\begin{array}{ccc}
 0& F_{3}^0 & -F_{2}^0\\
 -F_{3}^0& 0 & F_{1}^0\\
 F_{2}^0& -F_{1}^0 & 0
\end{array}\right),
\nonumber
\end{equation}

\begin{equation}
\bo B_4 =  \left(\begin{array}{ccc}
 F_2^0\kappa_2^0+F_3^0\kappa_3^0-B(d_{22}^0d_{11}^0-d_{21}^0d_{12}^0)& (F_3^0)'-F_1^0\kappa_2^0 & -(F_2^0)'-F_1^0\kappa_3^0\\
 -(F_3^0)'-F_2^0\kappa_1^0& F_3\kappa_3^0+F_1^0\kappa_1^0-B(d_{22}^0d_{11}^0-d_{21}^0d_{12}^0) & (F_1^0)'-F_2^0\kappa_3^0\\
 (F_2^0)'-F_3\kappa_1^0-B(d_{22}^0d_{31}^0-d_{21}^0d_{32}^0)& -(F_1^0)'-F_3^0\kappa_2^0+B(d_{12}^0d_{31}^0-d_{11}^0d_{32}^0) &
 F_1^0\kappa_1^0+F_2^0\kappa_2^0
\end{array}\right),
\nonumber
\end{equation}

\begin{equation}
\bo B_5 = \left(\begin{array}{ccc}
 d_{11}^0& d_{12}^0 & d_{13}^0\\
 d_{21}^0& d_{22}^0 & d_{23}^0\\
 d_{31}^0& d_{32}^0 & d_{33}^0
\end{array}\right),
\nonumber
\end{equation}

\begin{equation}
\bo B_6 = 2\omega \left(\begin{array}{ccc}
 d_{12}^0& -d_{11}^0 & 0\\
 d_{22}^0& -d_{21}^0 & 0\\
 d_{32}^0& -d_{31}^0 & 0
\end{array}\right).
\nonumber
\end{equation}

Matrices $\bo C_i$ appearing in equation (\ref{per:amomentum}) are given by

\begin{equation}
\bo C_1 = \bo B_1, \nonumber
\end{equation}

\begin{equation}
\bo C_2 = \left(\begin{array}{ccc}
 0& M_{3}^0 & -M_{2}^0\\
 -M_{3}^0& 0 & M_{1}^0\\
 M_{2}^0& -M_{1}^0 & 0
\end{array}\right),
\nonumber
\end{equation}

\begin{equation}
\bo C_3 = \left(\begin{array}{ccc}
 C_3^{11}& C_3^{12} & C_3^{13}\\
 C_3^{21}& C_3^{22} & C_3^{23}\\
 C_3^{31}& C_3^{32} & C_3^{33}
\end{array}\right),
\nonumber
\end{equation}
where

\begin{eqnarray}
C_3^{11}&=&M_3^0\kappa_3^0+M_2^0\kappa_2^0-P(\bo \omega \cdot \bo d_3^0)(\bo d_2^0 \times \bo \omega \cdot \bo d_1^0)-P(\bo \omega \cdot \bo d_2^0)(\bo d_3^0 \times \bo \omega \cdot \bo d_1^0), \nonumber \\
C_3^{12}&=&(M_3^0)'-M_1^0\kappa_2^0-PR(\bo \omega \cdot \bo d_1^0)(\bo d_1^0 \times \bo \omega \cdot \bo d_3^0), \nonumber\\
C_3^{13}&=&-(M_2^0)'-M_1^0\kappa_3^0-F_1^0+PR(\bo \omega \cdot \bo d_1^0)(\bo d_1^0 \times \bo \omega \cdot \bo d_2^0)+P(\bo \omega \cdot \bo d_1^0)(\bo d_2^0 \times \bo \omega \cdot \bo d_1^0), \nonumber\\
C_3^{21}&=&-(M_3^0)'-M_2^0\kappa_1^0+P(\bo \omega \cdot \bo d_2^0)(\bo d_2^0 \times \bo \omega \cdot \bo d_3^0), \nonumber\\
C_3^{22}&=&M_3^0\kappa_3^0+M_1^0\kappa_1^0+PR(\bo \omega \cdot \bo d_3^0)(\bo d_1^0 \times \bo \omega \cdot \bo d_2^0)+PR(\bo \omega \cdot \bo d_1^0)(\bo d_3^0 \times \bo \omega \cdot \bo d_2^0), \nonumber \\
C_3^{23}&=&(M_1^0)'-M_2^0\kappa_3^0-F_2^0-PR(\bo \omega \cdot \bo d_2^0)(\bo d_1^0 \times \bo \omega \cdot \bo d_2^0)-P(\bo \omega \cdot \bo d_2^0)(\bo d_2^0 \times \bo \omega \cdot \bo d_1^0), \nonumber\\
C_3^{31}&=&(M_2^0)'-M_3^0\kappa_1^0+F_1^0-P(\bo \omega \cdot \bo d_3^0)(\bo d_2^0 \times \bo \omega \cdot \bo d_3^0), \nonumber\\
C_3^{32}&=&-(M_1^0)'-M_3^0\kappa_2^0+F_2^0+PR(\bo \omega \cdot \bo d_3^0)(\bo d_1^0 \times \bo \omega \cdot \bo d_3^0), \nonumber\\
C_3^{33}&=&M_2^0\kappa_2^0+M_1^0\kappa_1^0+P(1-R)(\bo \omega \cdot \bo d_2^0)(\bo d_1^0 \times \bo \omega \cdot \bo d_3^0)+P(1-R)(\bo \omega \cdot \bo d_1^0)(\bo d_2^0 \times \bo \omega \cdot \bo d_3^0),\nonumber \\
\nonumber
\end{eqnarray}

\begin{equation}
\bo C_4 = \left(\begin{array}{ccc}
 0& -1 & 0\\
 1& 0 & 0\\
 0& 0 & 0
\end{array}\right),
\nonumber
\end{equation}

\begin{equation}
\bo C_5 = P\left(\begin{array}{ccc}
 1& 0 & 0\\
 0& R & 0\\
 0& 0 & (1+R)
\end{array}\right),
\nonumber
\end{equation}

\begin{equation}
\bo C_6 = 2P \left(\begin{array}{ccc}
 (\bo d_2^0 \times \omega \times \bo d_3^0) \cdot \bo d_1^0& 0 & -(\bo d_2^0 \times \omega \times \bo d_1^0) \cdot \bo d_1^0\\
 0& -R(\bo d_1^0 \times \omega \times \bo d_3^0) \cdot \bo d_2^0 & R(\bo d_1^0 \times \omega \times \bo d_2^0) \cdot \bo d_2^0\\
 (\bo d_2^0 \times \omega \times \bo d_3^0) \cdot \bo d_3^0& -R(\bo d_1^0 \times \omega \times \bo d_3^0) \cdot \bo d_3^0 & R(\bo d_1^0 \times \omega \times \bo d_2^0) \cdot \bo
 d_3^0-(\bo d_2^0 \times \omega \times \bo d_1^0) \cdot \bo d_3^0
\end{array}\right).
\nonumber
\end{equation}

Matrices $\bo D_i$ appearing in equation (\ref{per:constitutive}) are given by

\begin{equation}
\bo D_1 = \left(\begin{array}{ccc}
 -1& 0 & 0\\
 0& -R & 0\\
 0& 0 & -\frac{\Gamma(1+R)}{2}
\end{array}\right),
\nonumber
\end{equation}

\begin{equation}
\bo D_2 = \left(\begin{array}{ccc}
 0& \kappa_3^0 & -M_2^0- (1-R)\kappa_2^0\\
 - \kappa_3^0& 0 & M_1^0- (1-R)\kappa_1^0\\
 M_2^0- \kappa_2^0\left(R-\frac{\Gamma(1+R)}{2}\right)&
 -M_1^0+ \kappa_1^0\left(1-\frac{\Gamma(1+R)}{2}\right) &0
\end{array}\right),
\nonumber
\end{equation}

\begin{equation}
\bo D_3 = \gamma\left(\begin{array}{ccc}
 0& -\kappa_3^0 & \kappa_2^0\\
 R\kappa_3^0& 0 & -R\kappa_1^0\\
 -\frac{\Gamma(1+R)}{2}\kappa_2^0& \frac{\Gamma(1+R)}{2}\kappa_1^0 &0
\end{array}\right),
\nonumber
\end{equation}

\begin{equation}
\bo D_4 = -\gamma \bo D_1. \nonumber
\end{equation}

All the $\kappa_i^0$ in the above can be expressed in terms of the moments
$M_i^0$ by means of the constitutive relations (\ref{stat:constitutive}).

\newpage

\appendix
\begin{center}
\textbf{Appendix B: Helical solutions}
\end{center}

\noindent
If we {\it assume} a helical shape for the magnetically buckled rod then
we can derive exact solutions as well as buckling loads. For this it is
convenient to introduce Euler angles $(\theta$, $\psi$, $\phi)$
and relate the director frame $\{\bo{d}_1,\bo{d}_2,\bo{d}_3\}$ to the
rotating frame $\{\bo{e}_1,\bo{e}_2,\bo{e}_3\}$ as follows:
\begin{eqnarray}
\bo{d}_1&=&(\sin\phi\cos\psi+\cos\phi\cos\theta\sin\psi)\,\bo{e}_1 +
(\sin\phi\sin\psi-\cos\phi\cos\theta\cos\psi)\,\bo{e}_2 \nonumber \\
&&+ \cos\phi\sin\theta\,\bo{e}_3, \nonumber\\
\bo{d}_2&=&(\cos\phi\cos\psi-\sin\phi\cos\theta\sin\psi)\,\bo{e}_1 +
(\cos\phi\sin\psi+\sin\phi\cos\theta\cos\psi)\,\bo{e}_2 \nonumber \\
&&- \sin\phi\sin\theta\,\bo{e}_3, \label{angles} \\
\bo{d}_3&=&-\sin\theta\sin\psi\,\bo{e}_1 + \sin\theta\cos\psi\,\bo{e}_2
+ \cos\theta\,\bo{e}_3.\nonumber
\end{eqnarray}
The rate of change of the director frame is given by
\begin{equation}
{\bo{d}_i}'=\bo{\kappa}\times\bo{d}_i \quad\quad (i=1,2,3),
\label{strains}
\end{equation}
where $\bo{\kappa}=\kappa_1\bo{d}_1+\kappa_2\bo{d}_2+\kappa_3\bo{d}_3$
is the curvature vector. Inverting (\ref{strains}) and using orthogonality
of the directors gives
\begin{equation}
\kappa_1={\bo{d}_2}'\cdot\bo{d}_3, \quad\quad\quad
\kappa_2={\bo{d}_3}'\cdot\bo{d}_1, \quad\quad\quad
\kappa_3={\bo{d}_1}'\cdot\bo{d}_2,
\label{kappa}
\end{equation}
which, on inserting (\ref{angles}), yields
\begin{eqnarray}
&& \kappa_1=\psi'\sin\theta\cos\phi - \theta'\sin\phi, \nonumber \\
&& \kappa_2=-\psi'\sin\theta\sin\phi - \theta'\cos\phi, \label{curvs} \\
&& \kappa_3=\phi'+\psi'\cos\theta. \nonumber
\end{eqnarray}
So for the total curvature $\kappa$ we find
\begin{equation}
\kappa^2=\kappa_1^2+\kappa_2^2=\theta'^2+\psi'^2\sin^2\theta.
\label{curvature}
\end{equation}

Now we make the assumption of a helical centreline and uniformly rotating
directors by taking $\theta$ to be a constant and setting $\psi'=\Omega$,
$\phi'=\nu$ ($\Omega$ and $\nu$ constants). (That $\phi'$ may be taken
constant is not a priori clear. However, performing the following calculation
without this assumption, while still taking $\psi'$ constant, one is quickly
led to conclude that $\phi'$ is in fact constant, so to simplify the
presentation we assume constancy of $\phi'$ from the start.)
By integrating $\bo{x}'=\bo{d}_3$ we then find for the shape
\begin{equation}
\begin{array}{ll}
x(s)=r\cos(\Omega s+\psi(0)), \\
y(s)=r\sin(\Omega s+\psi(0)), \quad\quad \mbox{where} \quad\quad
\displaystyle r=\frac{\sin\theta}{\Omega}=\frac{\sin^2\theta}{\kappa}, \vspace*{0.2cm} \\
z(s)=L+(s-L)\cos\theta,
\end{array}
\label{shape}
\end{equation}
and with this solution the boundary conditions (\ref{eq:BC36}) --
(\ref{eq:BC47}) reduce to the four kinematical conditions
\begin{equation}
\cos(\Omega L+\psi(0))=0, \quad \cos(\chi+\psi(0))=0, \quad
\sin(\nu L+\phi(0))=0, \quad \sin\phi(0)=0,
\label{bc2}
\end{equation}
where we have allowed for a non-zero angle $\chi$ between the axes
$\bo v_0$ and $\bo v_1$ (see Fig.~\ref{fig:hel:BC}). Compatibility of the
director and fixed frame for the case of a straight (but twisted) rod (i.e.,
with $\theta=0$) requires
\begin{equation}
\psi(0)+\phi(0)=\frac{\pi}{2}-\chi, \quad\quad \psi(L)+\phi(L)=\frac{\pi}{2}
\label{bc3}
\end{equation}
(both frames are aligned at $s=L$ and have a relative rotation $\chi$ about
$\bo e_3=\bo k$ at $s=0$). Equations (\ref{bc2}) and (\ref{bc3}) imply
\begin{equation}
\Omega=\frac{n\pi+\chi}{L}, \quad\quad \nu=-\frac{n\pi}{L}, \quad\quad
\phi(0)=0,\pi, \quad\quad \psi(0)=\frac{\pi}{2}-\chi-\phi(0) \quad\quad
(n\in\mathbb{Z}).
\label{omega}
\end{equation}

To find the reduced equilibrium equation we integrate the force balance
equation (\ref{eq:lmomentum4_whirl}) to get (letting
$F_x=\bo{F}\cdot\bo{e}_1$, $F_y=\bo{F}\cdot\bo{e}_2$, $F_z=\bo{F}\cdot\bo{e}_3$)
\begin{equation}
\begin{array}{ll}
\displaystyle
F_x(s)=-IB_0 r\sin(\Omega s+\psi(0)) - \frac{r\rho\omega^2A}{\Omega}
\sin(\Omega s+\psi(0)), \vspace*{0.2cm} \\
\displaystyle
F_y(s)=IB_0 r\cos(\Omega s+\psi(0)) + \frac{r\rho\omega^2A}{\Omega}
\cos(\Omega s+\psi(0)), \vspace*{0.3cm} \\
\displaystyle
F_z(s)=T. \nonumber
\end{array}
\label{forces}
\end{equation}
and insert these expressions, together with (\ref{curvs}) and the
constitutive relations (\ref{eq:constitutive1}) with $\gamma_v=0$, in the
moment balance equation (\ref{eq:amomentum4_whirl}) written out in the
director frame. The twisting moment is directly given by its constitutive
relation to be
\begin{equation}
M_3=GJ\kappa_3=GJ(\phi'+\psi'\cos\theta)=GJ(\nu+\Omega\cos\theta),
\label{M_3}
\end{equation}
i.e., a constant. In the isotropic case ($I_1=I_2=:I_0$) the moment balance
equation for $M_3$ is therefore identically satisfied, while the equations
for $M_1$ and $M_2$ each give
\begin{equation}
(EI_0-GJ)\Omega\cos\theta-GJ\nu-\frac{IB_0}{\Omega^2}\cos\theta+\frac{T}{\Omega}-
\frac{\rho\omega^2A}{\Omega^3}\cos\theta+\frac{\rho\omega^2 I}{\Omega}\cos\theta=0,
\label{equil_eq}
\end{equation}
with $\Omega$ and $\nu$ given by (\ref{omega}). The torque $M$ of the lower
axis $\bo v_0$ about $\bo k$ is given by
\begin{equation}
M=M_z(0)=EI_0\Omega\sin^2\theta+GJ(\nu+\Omega\cos\theta)\cos\theta.
\end{equation}

Thus we have an infinite set of post-buckling helical solution branches
parametrised by $\theta\in [0,\pi/2]$. The bifurcating curves in
Fig.~\ref{fig:hel:bif-helix}, for instance, are obtained by plotting
(nondimensionalised) $r$ against $IB_0$ for $n=1,...,6$, using equations
(\ref{shape}) and (\ref{equil_eq}). The handedness of helical solutions is
determined by the sign of $\Omega$, right-handed for $\Omega>0$. The two
solutions for $\phi(0)$ in (\ref{omega}) correspond to the two branches
emanating from the pitchfork bifurcations, both with the same handedness
but with opposite signs for both $\kappa_1$ and $\kappa_2$.

The buckling condition is obtained by setting $\theta=0$:
\begin{equation}
IB_0=\frac{EI_0}{L^3}(n\pi+\chi)^3-\frac{GJ\chi}{L^3}(n\pi+\chi)^2+
\frac{T}{L}(n\pi+\chi)-\frac{\rho\omega^2AL}{n\pi+\chi}+\frac{\rho\omega^2I}{L}(n\pi+\chi),
\label{buckling}
\end{equation}
or, for $\chi=0$ and in dimensionless parameters,
\begin{equation}
B=n^3\pi^3+n\pi \bar{T}-\frac{\bar{\omega}^2}{n\pi}+n\pi\bar{\omega}^2P.
\label{B}
\end{equation}
This expression shows that an applied tension ($\bar{T}$) and inertia ($P$)
stiffen the rod against helical buckling, while whirl ($\bar{\omega}$)
softens it, as expected. The critical values of $B$ agree with the pitchfork
bifurcations in Figs~\ref{fig:hel:bif-helix} and \ref{fig:hel:bif-diag-w2},
and with (\ref{char_coat}) for the present case of an isotropic rod. Formula
(\ref{buckling}) is of course also valid in the absence of a magnetic field
($B_0=0$), in which case it gives critical values for buckling due to
compression ($T$), twist ($\chi$) or whirl ($\omega$).

In the anisotropic case ($I_1\neq I_2$) the moment balance equation for
$M_3$ becomes
$0=(E\Omega^2+\rho\omega^2)(I_2-I_1)\sin^2\theta\sin\phi\cos\phi$,
whose only solution compatible with the boundary conditions (\ref{omega})
is $\phi\equiv 0$ or $\pi$. But this implies $\nu=0$ and hence, by
(\ref{omega}), $n=0$ and $\Omega=\chi/L$. We conclude that in the case
$\chi=0$, i.e., in the case of parallel axes $\bo v_0$ and $\bo v_1$, an
anisotropic rod cannot buckle into
a helical solution. This agrees with our numerical results, which show that
the solutions bifurcating at the critical loads (\ref{char_coat}) when
$R\neq 1$ are non-helical.

\end{document}